\documentclass[12pt]{article}
\usepackage{latexsym}
\usepackage[all]{xy}
\usepackage{amsmath}
\usepackage{amssymb}
\usepackage{amsfonts}
\usepackage{color}
\usepackage{theorem}
\usepackage{mathtools}

\newcommand{\ql}{\mathbb{Q}_{\ell}}
\newcommand{\Ql}{\mathbb{Q}_{\ell}}
\newcommand{\rl}{r_{\ell}}
\newcommand{\Sing}{\mathsf{Sing}}
\newcommand {\Perf} {\mathsf{Perf}}
\newcommand {\Coh} {\mathsf{Coh}^{b}}
\newcommand{\Gal}{\mathrm{Gal}}
\newcommand{\Ind}{\mathsf{Ind}}
\newcommand{\Mod}{\mathsf{Mod}}

\newcommand {\Map} {\mathbb{R}\mathbf{Map}}

\newcommand {\rh} {\mathbb{R}\underline{Hom}}
\newcommand {\rch} {\mathbb{R}\underline{\mathcal{H}om}}
\newcommand {\OO} {\mathcal{O}}

\newcommand {\T} {\mathbb{T}}

\newcommand {\Spec} {\mathbf{Spec}}
\newcommand  {\dg}     {\mathbf{dg}}
\newcommand  {\dgcat}     {\mathbf{dgCat}}
\newcommand  {\dgCAT}     {\mathbf{dgCAT}}

\newcommand{\s}{\infty}

\newcommand{\BU}{\mathbf{B}\mathbb{U}}
\newcommand{\E}{\mathbb{E}}

\newcommand{\SH}{\mathcal{SH}}

\newcommand{\Alg}{\mathbf{Alg}}
\newcommand{\D}{\mathcal{D}}
\newcommand{\B}{\mathcal{B}}
\newcommand{\str}{\mathrm{str}}

\newtheorem{thm}{Theorem}[subsection]
\newtheorem{prop}[thm]{Proposition}
\newtheorem{lem}[thm]{Lemma}

\newtheorem{df}[thm]{Definition}
\newtheorem{cor}[thm]{Corollary}
\newtheorem{conj}[thm]{Conjecture}
\newtheorem{hyp}[thm]{Hypothesis}
{\theorembodyfont{\rmfamily} \newtheorem{rmk}[thm]{Remark}}

\begin{document}

\title{\textbf{Trace and K\"unneth formulas for singularity categories  and applications}}  

\author{Bertrand To\"en\footnote{Partially supported by ERC-2016-ADG-741501
and ANR-11-LABX-0040-CIMI within the program ANR-11-IDEX-0002-02.}\, and Gabriele 
Vezzosi}  

\date{\small{January 2019}}

\maketitle

\begin{abstract} \noindent We present an $\ell$-adic trace formula for
saturated and admissible dg-categories over a base monoidal dg-category. 
Moreover, we prove K\"unneth formulas for dg-category of singularities, and for inertia-invariant vanishing cycles.
As an application, we prove a version of
Bloch's Conductor Conjecture (stated by Spencer Bloch in 1985),  
under the additional hypothesis that the monodromy action of the inertia group is unipotent. \\

\end{abstract}

\tableofcontents

\section{Introduction}

Motivated by applications to arithmetic geometry, the first one being presented in this paper, we develop methods of derived and non-commutative geometry on a mixed-characteristic base. Though both derived and non-commutative geometry have been mainly used so far over a base field (possibly of positive characteristic), their flexibility is definitely wider, as we try to demonstrate in this paper. More precisely, we prove and investigate a trace formula for dg-categories, and K\"unneth formulas for dg-categories of singularities and for inertia-invariant vanishing cycles, and we give an application of these results to a categorical version of Bloch's Conductor Conjecture (BCC) in the presence of unipotent monodromy. We are convinced that the tools developed in this paper will have other applications to arithmetic geometry, some of which are already under investigation by the authors. Another seemingly more flexible tool, that of \emph{integration maps} from $K$-theory of dg-categories, will be developed in a forthcoming paper and applied to more general cases of the original BCC. \\

In this paper, we present four main results. As a first step, in Section 2, we prove a quite general \emph{trace formula for dg-categories} (or non-commutative schemes). This is done by using the non-commutative $\ell$-adic realization functor $r_{\ell}$ recently introduced in \cite{brtv} as a functor from dg-categories over an excellent discrete valuation ring $A$ to $\ell$-adic sheaves on $Spec \, A$. Via $r_{\ell}$, we first introduce a $\ell$-adic version of the Chern character for dg-categories over $A$, as a lax-monoidal natural transformation $Ch_{\ell}: \mathsf{HK} \to |r_{\ell}|$, where $\mathsf{HK}$ is the non-connective, homotopy invariant $K$-theory functor on dg-categories, and $|r_{\ell}|$ is the Eilenberg-Mac Lane construction applied to the derived global sections of the $\ell$-adic realization functor $r_{\ell}$.
We prove the following result (Theorem \ref{ttrace} in the paper).\\

\noindent \textbf{Theorem A (Trace formula for dg-categories).} \emph{Let $\B$ be a monoidal dg-category over $A$. For $T$ a dg-category which is a (left) $\B$-module, satisfying both a version of smoothness and properness \emph{over $\B$} and an admissibility property\footnote{The $\ell$-adic realization functor is only lax-monoidal, and we say that $T$ is $r_{\ell}$-admissible over $\B$ if $r_{\ell}$ behaves like a symmetric monoidal functor for $T$ over $\B$ (see Definition \ref{ammisibilta}).} with respect to $r_{\ell}$,  a trace formula holds \begin{equation}\label{traceintro}Ch_{\ell}([\mathsf{HH}(T/\B; f)]) = tr_{r_{\ell}(\B)}(r_{\ell}(f))\end{equation}
for any endomorphism $f:T \to T $ of $\B$-modules.}\\

\noindent   Here $[\mathsf{HH}(T/\B; f)]$ is the class in $\mathsf{HK}_0(\mathsf{HH}(\B/A))$ induced by the endomorphism $f$, and $\mathsf{HK}_0(\mathsf{HH}(\B/A))$ is the $0$-th $K$-theory of Hochschild homology of $\B$ over $A$. The r.h.s.  $tr_{r_{\ell}(\B)}(r_{\ell}(f))$ is the trace\footnote{Or rather its image under the canonical map $\alpha: H^0(S_{\textrm{\'et}}, \mathsf{HH}(r_{\ell}(\B)/r_{\ell}(A))) \to H^0(S_{\textrm{\'et}}, r_{\ell}(\mathsf{HH}(\B/A)))$.} of the endomorphism $r_{\ell}(f)$, and the trace formula (\ref{traceintro}) is an equality in 
$\pi_0 (|\rl(\mathsf{HH}(B/ A))|) \simeq H^0(S_{\textrm{\'et}}, \rl(\mathsf{HH}(B/ A)) )$. \\ 
Note that the extra generality of working over $\B$ (instead of over $A$, or over an $E_{\infty}$- algebra over $A$) is actually needed for applications, as shown in Section 5: the $\ell$-adic realization of a dg-category is always $2$-periodic, i.e. it is a module over $\Ql(\beta):= \oplus_{n\in \mathbb{Z}}\Ql(2n)$, so it has no chance of being smooth and proper over $A$ itself, unless it is trivial; in particular, no reasonable trace formula is available over $A$.\\

Our second main result (Section 3) is a \emph{K\"unneth type formula for inertia-invariant vanishing cycles}. Here we push forward the investigation begun in \cite{brtv} about the relation between vanishing cycles and the $\ell$-adic realization of the dg-category of singularities. For $X$ and $Y$ regular schemes, endowed with a flat, proper and generically smooth map to the strictly henselian excellent trait $S= Spec \, A$, we define an inertia-invariant convolution product $(E\circledast F)^\mathrm{I}$, for $E$ an $\ell$-adic sheaf on the special fiber $X_s$, and $F$ an $\ell$-adic sheaf on the special fiber $Y_s$ ($\mathrm{I}:= \mathrm{Gal}(\bar{K}/K)$ denoting the inertia group in this situation). If we denote by $ \nu_X$ (respectively, $\nu_Y$) the complex of vanishing cycles for $X/S$ (respectively, $Y/S$), we prove the following result  (Theorem \ref{kunn} in the paper).\\

\noindent \textbf{Theorem B (K\"unneth formula for inertia-invariant vanishing cycles).} \emph{We have equivalences $$(\nu_X \circledast \nu_Y)^{\mathrm{I}} \simeq r_{\ell}(\Sing(X\times_S Y)) \simeq \textrm{Cofib}(\eta_{X\times_S Y} :\Ql(\beta) \longrightarrow \omega_{X\times_S Y}(\beta)).$$ where $\Sing(X\times_S Y)=\Coh(X\times_S Y)/\Perf (X\times_S Y)$ is the dg-category of singularities of the fiber product $X\times_S Y$, and $\eta_{X\times_S Y}$ is the (2-periodized) $\ell$-adic fundamental class.} \\

This is our K\"unneth formula for inertia-invariant vanishing cycles, and it seems to be a new result in the theory of vanishing cycles, especially in the mixed characteristic case. Note that it might also be viewed as a Thom-Sebastiani formula for inertia-invariant vanishing cycles but we would like to stress that it is not a consequence of the usual Thom-Sebastiani formula for vanishing cycles (see \cite{Ilts}), and that it holds only if inertia invariants are taken into account in defining the convolution $(E\circledast F)^\mathrm{I}$\footnote{In particular, only $(E\circledast F)^\mathrm{I}$, and not $(E\circledast F)$, makes sense in our context.}. We also discuss appropriate conditions  ensuring that $(\nu_X \circledast \nu_Y)^{\mathrm{I}}$ is equivalent to $(\nu_X \boxtimes \nu_Y)^{\mathrm{I}}$ (Corollary \ref{ckunn}). \\

The third main result of this paper is a \emph{K\"unneth type formula for dg-categories of singularities} (Section 4). For $X$ and $Y$ regular schemes, endowed with a flat, proper and generically smooth map to the strictly henselian excellent trait $S= Spec \, A$, we may consider the dg-categories of singularities $\Sing(X_s):= \Coh(X_s)/\Perf(X_s)$ and $\Sing(Y_s):=\Coh(X_s)/\Perf(X_s)$ of the corresponding special fibers. Consider $\B:=\Sing(s\times_S s)$ ($s$ being the closed point in $S$, and $s \times_S s$ being the derived fiber product); then $\B$ is a monoidal dg-category for the convolution product coming from the derived groupoid structure of $s \times_S s$. Moreover, $\B$ acts on both $\Sing(X_s)$ and $\Sing(Y_s)$, in such a way that the tensor product $\Sing(X_s)^{o}\otimes_{\B} \Sing(Y_s)$ makes sense as a dg-category over $A$ (Proposition \ref{cotens}). Our K\"unneth formula for dg-categories of singularities is then the following result (Theorem \ref{dg-kunn} in the paper).\\

\noindent \textbf{Theorem C (K\"unneth formula for dg-categories of singularities).} \emph{There is a canonical equivalence  $$\Sing(X_s)^{o}\otimes_{\B} \Sing(Y_s) \simeq \Sing(X\times_S Y)$$ as dg-categories over $A$.}\\ 

\noindent Note that this result is peculiar to singularity categories: it is false if we replace $\Sing$ by $\Coh$. Also notice that $\B$ is defined as the convolution dg-category of a derived groupoid, and is an object in derived algebraic geometry that \emph{cannot} be described within classical algebraic
geometry. Finally, we prove that $\Sing (X_s)$ is smooth and proper (i.e. saturated) over $\B$, for any $X$ regular scheme, endowed with a flat, proper and generically smooth map to $S$.  This fact is well-known in characteristic zero (see, for instance, \cite{pre}) but is, in our opinion, 
a deep and surprising result in our general setting. \\Smoothness and properness over $\B$ is one half of the properties needed to apply the trace formula to $\Sing (X_s)$ over $\B$. Note that, without further hypothesis, it is not true that $\Sing (X_s)$ is also admissible (with respect to $r_{\ell}$).\\

Our fourth and final main result is a \emph{categorical version of Bloch's Conductor formula for unipotent monodromy} (section 5).\\
In his seminal 1987 paper \cite{bl}, S. Bloch introduced 
what is now called Bloch's intersection number $[\Delta_X,\Delta_X]_S$, for
a flat, proper map of schemes $X \longrightarrow S$ where $X$ is regular, and $S$ is a henselian trait $S$. This number
can be defined as the degree of the localized top Chern class of the coherent sheaf 
$\Omega^1_{X/S}$ and measures the ``relative'' singularities of $X$ over $S$. 
In the same paper, Bloch introduced his famous \emph{conductor formula}, which 
can be seen as a conjectural computation of the Bloch's intersection number in terms
of the arithmetic geometry of $X/S$. It reads as follows.\\

\noindent \textbf{Bloch's Conductor Conjecture (BCC)}\\
\emph{We have an equality
$$[\Delta_X,\Delta_X]_S = \chi(X_{\bar{k}}) - 
\chi(X_{\bar{K}}) - Sw(X_{\bar{K}}),$$
where $X_{\bar{k}}$ and $X_{\bar{K}}$ denotes the special 
and generic geometric fibers of $X$ over $S$, $\chi(-)$ denotes
$\Ql$-adic Euler characteristic, and $Sw(-)$ is the Swan conductor.}\\

In \cite{bl} the above formula is proven in relative dimension $1$. Further
results implying special cases of BCC have been obtained 
since then by Kato-Saito \cite{ks} and others, 
the most recent one being a full proof in the geometric case by T. Saito \cite{sait} (see section \S 5 for a more detailed discussion about the status of the art about BCC).
In the mixed characteristic case, the conjecture is open in general
outside the cases covered in \cite{ks}. In particular, for isolated
singularities the above conjecture already appeared in Deligne's expos\'e
\cite[Exp. XVI]{sga7}, and remains open.  \\
In Section 5 we propose a first step 
towards a new understanding of Bloch's conductor formula using the results developed
in the previous sections.
We start by defining an analog of Bloch's intersection
number, that we call the \emph{categorical Bloch's intersection
class} (Definition \ref{bnum}), and denote it by $[\Delta_X,\Delta_X]_S^{\textrm{cat}}$.
It is an element in $H^0(S_{\textrm{\'et}}, r_{\ell}(\mathsf{HH}(\B/A)))$, where $\B=\Sing(s\times_S s)$ (as in Section 4, see above), and 
it is basically defined as an intersection class in the setting of non-commutative
algebraic geometry. The precise comparison 
with the original Bloch's number is not covered in this work and
will appear in a forthcoming paper. \\ It is easy to see that there are canonical inclusions $\mathbb{Z}\hookrightarrow \Ql \hookrightarrow H^0(S_{\textrm{\'et}}, \mathsf{HH}(r_{\ell}(\B)/r_{\ell}(A)))$, and for any $\lambda \in \mathbb{Z}$, we will write $\lambda^{\wedge}$ its image under the canonical map $\alpha: H^0(S_{\textrm{\'et}}, \mathsf{HH}(r_{\ell}(\B)/r_{\ell}(A))) \to H^0(S_{\textrm{\'et}}, r_{\ell}(\mathsf{HH}(\B/A)))$.

Our fourth main result is then the following theorem (Theorem \ref{blochthm} in the paper).\\

\noindent \textbf{Theorem D (Categorical BCC for unipotent monodromy)} \emph{With the notations of BCC, if the monodromy action of the inertia on $H^*(X_{\bar{K}},\Ql)$
is unipotent, then we have an equality 
$$[\Delta_X,\Delta_X]^{\textrm{cat}}_S = \chi(X_{\bar{k}})^{\wedge} - 
\chi(X_{\bar{K}})^{\wedge}.$$}

\noindent Since unipotent monodromy action implies tame action, the Swan conductor vanishes for unipotent monodromy, so that  Theorem D is completely analogous to BCC under this hypothesis. \\
We believe the main ideas in our proof of Theorem D are new in the subject, and might be also useful  
to answer other related questions in algebraic and arithmetic geometry. 
The key point in the proof of Theorem D is that it is a \emph{direct 
consequence} of the trace formula for dg-categories (Theorem A). 
By the main result of \cite{brtv}, the $\ell$-adic realization  of the dg-category of singularities $\Sing (X_s)$ of the special fiber is the inertia invariant part of vanishing cohomology of $X/S$ (suitably 2-periodized). Moreover, 
our K\"unneth formulas for inertia-invariant vanishing cycles (Section 3), and for dg-categories of singularities (Section 4), imply that $\Sing (X_s)$ is smooth and proper over $\B$. Finally, again our K\"unneth formulas \emph{together with the hypothesis of unipotent monodromy} show that $\Sing (X_s)$ is also $r_{\ell}$-admissible, so that we may apply our trace formula to the identity endomorphism of $\Sing (X_s)$ over $\B$: Theorem D is then 
an immediate corollary of the trace formula applied to the identity map of $\Sing (X_s)$.\\

We conclude this introduction by remarking that the hypothesis of unipotent inertia action in Theorem D is, of course, a bit restrictive: unipotent monodromy implies tame monodromy and our theorem does not deal with the interesting arithmetic
aspects encoded by the Swan conductor. However we
are convinced that we can go beyond the unipotent case, and prove some new cases of BCC conjecture (possibly the whole unrestricted  
conjecture), by considering $\Sing(X\times_S S')$ as a module over $\B':= \Sing(S'\times_S 
S')$, where $S'=Spec \, A'\to S$ is a totally ramified Galois extension of strictly henselian traits such that the inertia for $S'$ acts unipotently on $H^*((X\times_S S')_{\bar{K'}},\Ql)$. Here $K'$ is the fraction field of $A'$, and the existence of such an extension $S'\to S$ is guaranteed by Grothendieck local monodromy theorem (\cite[Exp. I, Th\'eor\`eme 1.2]{sga7}).
This strategy for BCC, as well as the comparison between the categorical and the classical
Bloch numbers, will be investigated in a future work. \\

\noindent \textbf{Acknowledgments.} We thank A. Abb\`es, A. Blanc, M. Robalo and T. Saito, for various helpful conversations 
on the topics of this paper.
We also thank all the french spiders and the south-italy sea waves in the 2017 summer nights, for keeping the authors company 
 while they were, separately, writing down a first version of this manuscript.

This work is
part of a project that has received funding from the European Research Council (ERC) 
under the 
European
Union's Horizon 2020 Research and Innovation Programme (Grant Agreement No 741501). \\

\bigskip

\noindent \textbf{Notations.} \begin{itemize}
\item Throughout the text, $A$ will denote a (discrete) commutative noetherian ring. 
When needed, $A$ will be required to satisfy further properties (such as being 
excellent, local and henselian) that will be made precise in due course.
\item $L(A)$ will denote the $A$-linear dg-category of (fibrant and) cofibrant $A$-dg-modules  
localized with respect to quasi-isomorphisms.
\item $\dgcat_A$ will denote the Morita $\s$-category of dg-categories over $A$ 
(see \S 2.1). 
\item $\mathbf{Top}$ denotes the $\infty$-category of spaces (obtained, e.g. as the 
coherent nerve of the Dwyer-Kan localization of the category of simplicial sets along 
weak homotopy equivalences). $\mathbf{Sp}$ denotes the $\infty$-category of spectra.
\item If $R$ is a an associative and unital monoid in a symmetric monoidal $\infty$-category $\mathcal{C}$,
we will write $\mathbf{Mod}(R)$ or $\mathbf{Mod}_{\mathcal{C}}(R)$ for the $\infty$-category of left $R$-modules in $\mathcal{C}$ (\cite{luha}).

\end{itemize}

\section{A non-commutative trace formalism}

\subsection{$\s$-Categories of dg-categories} 

We denote by $A$ a commutative ring. We remind here some basic
facts about the $\s$-category of dg-categories, its monoidal structure
and its theory of monoids and modules.  \\

We consider the category $dgCat_{A}$ of small $A$-linear dg-categories and $A$-linear dg-functors.
We remind that an $A$-linear dg-functor $T \longrightarrow T'$ is a Morita equivalence
if  the induced functor of the corresponding derived categories of dg-modules
$f^* : D(T') \longrightarrow D(T)$
is an equivalence of categories (see \cite{to-dgcat} for details). The $\s$-category of dg-categories over $S$
is defined to be the localisation of $dgCat_{A}$ along these Morita equivalences, 
and will be denoted by $\dgcat_S$ or $\dgcat_A$. Being the $\s$-category associated
to a model category, $\dgcat_A$ is a presentable $\s$-category.
As in \cite[\S \, 4]{to-dgcat}, the tensor product of $A$-linear dg-categories can be derived to 
a symmetric monoidal structure on the $\s$-category $\dgcat_A$. This symmetric monoidal
structure moreover distributes over colimits making $\dgcat_A$ into 
a presentable symmetric monoidal $\s$-category.
We have a notion of rigid, or, equivalently, dualizable, object in 
$\dgcat_A$. It is a well known fact that 
dualizable objects in $\dgcat_A$ are precisely smooth and proper
dg-categories over $A$ (see \cite[Prop. 2.5]{to2}).

The compact objects in $\dgcat_A$ are the dg-categories of finite type over $A$ in the
sense of \cite{tv}. We denote their full sub-$\s$-category by $\dgcat_A^{ft} \subset \dgcat_A$.
The full sub-category $\dgcat_S^{ft}$ is preserved by  
the monoidal structure, and moreover any dg-category is a filtered colimit of 
dg-categories of finite type. 
We thus have a natural equivalence of symmetric monoidal $\s$-categories
$$\dgcat_A \simeq \mathbf{Ind}(\dgcat_S^{ft}).$$

We will from time to time have to work in a bigger $\s$-category, denoted
by $\dgCAT_A$, that contains $\dgcat_A$ as a non-full
sub-$\s$-category. By \cite{to2}, we have a symmetric monoidal $\s$-category 
$\dgcat_A^{lp}$ of \emph{presentable} dg-categories over $A$. We define
$\dgCAT_A$ the full sub-$\s$-category of $\dgcat_A^{lp}$
consisting of all compactly generated dg-categories. 
The $\s$-category $\dgcat_A$ can be identified with the non-full sub-$\s$-category
of $\dgCAT_A$ which consists of compact objects preserving dg-functors. This
provides a faithful embedding of symmetric monoidal $\s$-categories
$$\dgcat_A \hookrightarrow \dgCAT_A.$$
At the level of objects, this embedding sends a small dg-category 
$T$ to the compactly generated dg-category $\hat{T}$ of dg-modules over 
$T^{o}$. An equivalent description of $\dgCAT_A$ is as
the $\s$-category of small dg-categories together with the mapping spaces given 
by the classifying space of \emph{all} bi-dg-modules between small dg-categories. \\

\begin{df}\label{weakmon}
A \emph{monoidal $A$-dg-category} is a unital and associative monoid in the
symmetric monoidal $\s$-category $\dgcat_A$. A \emph{module over a monoidal $A$-dg-category $B$} 
will, by definition, mean a left $B$-module in 
$\dgcat_A$ in the sense of \cite{luha}, and the $\s$-category of (left) $B$-modules will be denoted
by $\dgcat_B$. 
\end{df}

For a $B$-module $T$, we have a morphism $\mu: B\otimes_A T \to T $ in 
$\dgcat_A$, that will be simply denoted by $(b,x) \mapsto b\otimes x$.
For a monoidal $A$-dg-category $B$, we will denote by $B^{\otimes\textrm{-op}}$ the
monoidal $A$-dg-category where the monoid structure is the opposite to the one of $B$, i.e. 
$b\otimes^{op} b' := b'\otimes b$. Note that $B^{\otimes\textrm{-op}}$ should not be confused with 
$B^o$ (which is still a monoidal $A$-dg-category), where the ``arrows'' and not the monoid 
structure have been reversed, i.e. $B^o(b,b'):= B(b', b)$. By definition, a \emph{right} $B$-module 
is a (left) $B^{\otimes\textrm{-op}}$-module. The $\s$-category of right $B$-modules
will be denoted by 
$\dgcat_{B^{\otimes\textrm{-op}}}$, or simply by $\dgcat^B$.
If $B$ is a monoidal $A$-dg-category, then $B^{\otimes\textrm{-op}}\otimes_A B$ is again a 
monoidal $A$-dg-category, and $B$ can be considered either as a left $B^{\otimes\textrm{-op}}\otimes_A 
B$ (denoted by $B^L$), or as a right $B^{\otimes\textrm{-op}}\otimes_A B$-module (denoted by $B^R$).
For $T$ a $B$-module, and $T'$ a right $B$-module, then $T'\otimes_A T$ is naturally a right 
$B^{\otimes\textrm{-op}}\otimes_A B$-module, and we define $$T'\otimes_B T := (T'\otimes_A T) 
\otimes_{B^{\otimes\textrm{-op}}\otimes_A B}B^L  $$ which is an object in $\dgcat_A$. \\

Let $B$ be a monoidal $A$-dg-category. We can consider the symmetric monoidal
embedding $\dgcat_A \hookrightarrow \dgCAT_A$, so that the image $\widehat{B}$ of $B$ 
is a monoid in $\dgCAT_A$. The $\s$-category of $\widehat{B}$-modules in 
$\dgCAT_A$ is denoted by 
$\dgCAT_B$, and its objects are called \emph{big $B$-modules}. The natural
$\s$-functor $\dgcat_B \longrightarrow \dgCAT_B$ is faithful and 
its image consists of all big $B$-modules $\widehat{T}$ such that the morphism
$\widehat{B}\hat{\otimes}_A \widehat{T} \longrightarrow \widehat{T}$ is a small morphism. 

It is known (see \cite{to2}) that the symmetric monoidal $\s$-category $\dgCAT_A$
is rigid, that for any $\widehat{T}$ its dual is 
given by $\widehat{T^o}$, and that the evaluation and coevaluation
morphisms are defined by $T$ considered as $T^o\otimes_A T$-module. 
This formally 
implies that if $\widehat{T}$ is a big $B$-module, then its dual
$\widehat{T^o}$ is naturally a \emph{right} big $B$-module. We thus have 
two big morphisms
$$\mu : \widehat{B} \hat{\otimes}_A \widehat{T} \longrightarrow \widehat{T} \,\, , \quad 
\mu^o : \widehat{T^o} \hat{\otimes}_A \widehat{B} \longrightarrow \widehat{T^o}.$$
These morphisms also provide a third big morphism
$$h : \widehat{T^o} \widehat{\otimes}_A \widehat{T} \longrightarrow \widehat{B}.$$
The big morphism $h$ is obtained by duality from 
$$\mu^* : \widehat{T} \longrightarrow 
\widehat{B} \hat{\otimes}_A \widehat{T}$$
the right adjoint to $\mu$.
We now make the following definitions.

\begin{df}\label{dcotens}
Let $B$ be a monoidal dg-category, and $T$ a $B$-module. We say that: 
\begin{enumerate}

\item $T$ is \emph{cotensored} (over $B$) if the big morphism
$\mu^o$ defined above is a small morphism (i.e. a morphism in $\dgcat_A$);

\item $T$ is \emph{proper} or \emph{enriched} (over $B$) if the big morphism
$h$ defined above is a small morphism (i.e. a morphism in $\dgcat_A$).

\end{enumerate}
\end{df}

We can make the above definition more explicit as follows. 
Let $B$ and $T$ be as above; for two objects $b \in B$ and $x \in T$, 
we can consider the dg-functor
$$x^b : T^o \longrightarrow L(A)$$
sending $y \in T$ to $T(b\otimes y),x)$ where  
$(b, x)\mapsto b\otimes x \, : B \otimes_A T \longrightarrow T$ is the $B$-module structure on $T$.
Then, $T$ is cotensored over $B$ if and only if for all $b$ and $x$ 
the above dg-module $x^b: T^o \longrightarrow L(A)$ is compact in the derived
category $D(T^o)$ of all $T^o$-dg-modules. When $T$ is furthermore assumed to be triangulated,
then this is equivalent to ask for the dg-module to be representable by 
an object $x^b \in T$.
In a similar manner, we can phrase properness of $T$ over $B$ by 
saying that, for any $x\in T$ and $y\in T$, the dg-module
$$B^o \longrightarrow L(A)$$
sending $b$ to $T(b\otimes x,y)$ is compact (or representable, if $B$
is furthermore assumed to be triangulated). \\

\begin{rmk} It is important to notice that, by definition, when $T$ is cotensored, then the big right $B$-module 
$\widehat{T^o}$ is in fact a small right $B$-module. Equivalently, 
when $T$ is cotensored, then $T^o$ is naturally a right $B$-module
inside $\dgcat_A$, with the right module structure given by the morphism in $\dgcat_A$
$$\mu^o : T^o\otimes_A B \longrightarrow T^o$$ sending $(x,b) \in T^o\otimes_A B$ to 
the cotensor $x^b \in T^o$.
\end{rmk}

We end this section by recalling some facts about existence of tensor
products of modules over monoidal dg-categories. As a general fact, since 
$\dgcat_A$ is a presentable symmetric monoidal $\s$-category, 
for any monoidal dg-category $B$ there exists a tensor product $\s$-functor
$$\otimes_B : \dgcat_B \times \dgcat^B \longrightarrow \dgcat_A,$$
sending a left $B$-module $T$ and a right $B$-module $T'$ 
to $T'\otimes_B T$ (see \cite{luha}). Now, if  
$T$ is a $B$-module which is also cotensored (over $B$) in the sense
of Definition \ref{dcotens}, we have that $T^o$ is a right $B$-module, and we can thus form
$$T^o \otimes_B T \in \dgcat_A.$$
When $T$ is not cotensored, the object $T \otimes_B T^o$ does not 
exist anymore. However, we can always consider the
presentable dg-categories $\widehat{T}$ and $\widehat{T^o}$
as left and right modules over $\widehat{B}$, respectively, and their tensor product
$\widehat{T^o} \widehat{\otimes}_{\widehat{B}}\widehat{T}$ now only makes sense
as a presentable dg-category which has no reason to be compactly generated, in general.
Of course, when $T$ is cotensored, this presentable dg-category is 
compactly generated and we have
$$\widehat{T^o \otimes_B  T} 
\simeq \widehat{T^o} \widehat{\otimes}_{\widehat{B}} \widehat{T}.$$

\begin{rmk}\label{little} The following, easy observation, will be useful in the sequel. Let $B$ be 
a monoidal dg-category, and assume that $B$ is generated, as a triangulated
dg-category, by its unit object $1 \in B$. Then, any big $B$-module is
small (i.e. in the image of $\dgcat_B \longrightarrow \dgCAT_B$), and also cotensored. 
\end{rmk}

\subsection{The $\ell$-adic realization of dg-categories}

We denote by $\SH_S$ the stable $\mathbb{A}^1$-homotopy $\s$-category of schemes over $S$ (see \cite[Def. 5.7]{vo} and \cite[\S\, 2]{ro}). 
It is a  presentable symmetric monoidal $\s$-category whose monoidal structure will be denoted by 
$\wedge_S$. Homotopy invariant algebraic K-theory defines an $\E_{\infty}$-ring 
object  in $\SH_S$ that we denote by $\BU_S$ (a more standard notation is $KGL_S$). We denote by 
$\mathbf{Mod}_{\SH_S}(\BU_S)$ the $\s$-category of $\BU_S$-modules in $\SH_S$. It is a presentable symmetric monoidal
$\s$-category whose monoidal structure will be denoted by $\wedge_{\BU_S}$. 

As proved in \cite{brtv}, there exists a lax symmetric monoidal $\s$-functor
$$\mathsf{M}^{-}: \dgcat_S \longrightarrow \mathbf{Mod}_{\SH_S}(\BU_S),$$
which will be denoted by $T \mapsto \mathsf{M}^T$ (while it is denoted by $T \mapsto \mathcal{M}_S^{\vee}(T)$ in \cite{brtv}). The precise construction of the $\s$-functor $\mathsf{M}^{-}$
is rather involved and uses in an essential manner the theory of non-commutative motives 
of \cite{ro} as well as the comparison with the stable homotopy theory of schemes. Intuitively, 
the $\s$-functor $\mathsf{M}^{-}$ sends a dg-category 
 $T$ to the homotopy invariant K-theory functor $S' \mapsto \mathsf{HK}(S' \otimes_S T)$. 
To be more precise, there is an obvious forgetful $\s$-functor 
$$\mathsf{U}: \mathbf{Mod}_{\SH_S}(\BU_S) \longrightarrow \mathbf{Fun}^{\s}(Sm_S^{op},\mathbf{Sp}),$$
to the $\s$-category of presheaves of spectra on the category $Sm_S$ of smooth $S$-schemes. For a given 
dg-category $T$ over $S$, the presheaf $\mathsf{U}(\mathsf{M}^T)$ is defined 
 by sending a smooth $S$-scheme $S'=\Spec\, A' \rightarrow \Spec\, A=S$
to $\mathsf{HK}(A'\otimes_A T)$, the homotopy invariant non-connective K-theory spectrum
of $A'\otimes_A T$ (see \cite[4.2.3]{ro}). \\

The $\s$-functor $\mathsf{M}^{-}$ satisfies some basic properties which we recall here.
\begin{enumerate}
\item The $\s$-functor $\mathsf{M}^{-}$ is a localizing invariant, i.e. for any 
short exact sequence $T_0 \hookrightarrow T \longrightarrow T/T_0$ of dg-categories over $A$, the induced sequence 
$$\xymatrix{\mathsf{M}^{T_0} \ar[r] & \mathsf{M}^T \ar[r] & \mathsf{M}^{T/T_0}}$$
exhibits $\mathsf{M}^{T_0}$ has the fiber of the morphism $\mathsf{M}^T \rightarrow \mathsf{M}^{T/T_0}$ in $\mathbf{Mod}(\BU_S)$.

\item The natural morphism
$\BU_S \longrightarrow \mathsf{M}^{A},$
induced by the lax monoidal structure of $\mathsf{M}^{-}$, is an equivalence of $\BU_S$-modules. 

\item The $\s$-functor $T \mapsto \mathsf{M}^T$ commutes with filtered colimits.

\item For any quasi-compact and quasi-separated scheme $X$, and any morphism  $p : X \longrightarrow S$, 
we have a natural equivalence of $\BU_S$-modules
$$\mathsf{M}^{\Perf(X)} \simeq p_*(\BU_X),$$
where $p_* : \mathbf{Mod}_{\SH_X}(\BU_X) \longrightarrow \mathbf{Mod}_{\SH_S}(\BU_S)$ is the direct image of
$\BU$-modules, and $\Perf(X)$ is the dg-category of perfect complexes on $X$.

\end{enumerate}

We now let $\ell$ be a prime number invertible in $A$. We denote by 
$Sh^{\mathsf{ct}}_{\Ql}(S)$ the $\s$-category of \emph{constructible} 
$\mathbb{Q}_{\ell}$-adic complexes on the \'etale site $S_{\acute{e}t}$ of $S$. It is a symmetric monoidal $\s$-category, and
we denote by 
$$Sh_{\Ql}(S):=\mathbf{Ind}(Sh^{\mathsf{ct}}_{\Ql}(S))$$
its completion under filtered colimits (see \cite[Def. 4.3.26]{gl}). 
According to \cite[Cor. 2.3.9]{ro}, there exists an $\ell$-adic realization $\s$-functor
$$R_{\ell} : \SH_S \longrightarrow Sh_{\Ql}(S).$$
By construction, $R_{\ell}$ is a symmetric monoidal $\s$-functor 
sending a smooth scheme $p : X \longrightarrow S$ to $p_!p^{!}(\mathbb{Q}_{\ell})$, 
or, in other words, to the relative $\ell$-adic homology of $X$ over $S$. 

We let $\T:=\mathbb{Q}_{\ell}[2](1)$, and we consider the $\E_{\infty}$-ring 
object in $Sh_{\Ql}(S)$
$$\Ql(\beta):=\oplus_{n\in \mathbb{Z}}\T^{\otimes n}.$$
 In this notation, $\beta$ stands for $\T$, 
and $\Ql(\beta)$ for the algebra of Laurent polynomials in $\beta$, so we might as well write
$$\Ql(\beta)=\Ql[\beta,\beta^{-1}].$$
As shown in \cite{brtv}, there exists a canonical equivalence $R_{\ell}(\BU_S) \simeq \Ql(\beta)$ of $\E_{\s}$-ring objects
in $Sh_{\Ql}(S)$,
that is induced by the Chern character from algebraic K-theory to 
\'etale cohomology. We thus obtain a well-defined symmetric monoidal $\s$-functor
$$\mathsf{R}_{\ell} : \mathbf{Mod}_{\SH_S}(\BU_S) \longrightarrow \mathbf{Mod}_{Sh_{\Ql}(S)}(\Ql(\beta)),$$
from $\BU_S$-modules in $\SH_S$ to $\Ql(\beta)$-modules in $Sh_{\Ql}(S)$.
By pre-composing with the functor $T \mapsto \mathsf{M}^T$, we obtain
a lax monoidal $\s$-functor 
$$\rl := \mathsf{R}_{\ell} \circ \mathsf{M}^{-}: \dgcat_S \longrightarrow \mathbf{Mod}_{Sh_{\Ql}(S)}(\Ql(\beta)).$$

\begin{df}\label{d1}
The $\s$-functor defined above 
$$\rl : \dgcat_S \longrightarrow \mathbf{Mod}_{Sh_{\Ql}(S)}(\Ql(\beta))$$
is called the \emph{$\ell$-adic realization functor for dg-categories over $S$}.
\end{df}  

From the standard properties of the functor $T \mapsto \mathsf{M}^T$, recalled above, we obtain the following
properties for the $\ell$-adic realization functor $T \mapsto \rl(T)$. 
 
\begin{enumerate}
\item The $\s$-functor $\rl$ is a \emph{localizing invariant}, i.e. for any 
short exact sequence $T_0 \hookrightarrow T \longrightarrow T/T_0$ of dg-categories over $A$, the induced sequence 
$$\xymatrix{\rl(T_0) \ar[r] & \rl(T) \ar[r] & \rl(T/T_0)}$$
is a fibration sequence in $\mathbf{Mod}(\Ql(\beta))$.

\item The natural morphism
$$\Ql(\beta) \longrightarrow \rl(A),$$
induced by the lax monoidal structure, is an equivalence in $\mathbf{Mod}_{Sh_{\Ql}(S)}(\Ql(\beta))$. 

\item The $\s$-functor $\rl$ commutes with filtered colimits.

\item For any separated morphism of finite type  $p : X \longrightarrow S$, 
we have a natural morphism of $\Ql(\beta)$-modules
$$\rl(\Perf(X)) \longrightarrow p_*(\Ql(\beta)),$$
where $p_* : \mathbf{Mod}_{Sh_{\Ql}(S)}(\mathbb{Q}_{\ell, X}(\beta)) \longrightarrow \mathbf{Mod}_{Sh_{\Ql}(S)}(\Ql(\beta))$ is induced by the
direct image $Sh^{\mathsf{ct}}_{\Ql}(X) \longrightarrow Sh^{\mathsf{ct}}_{\Ql}(S)$ of 
constructible $\Ql$-complexes.
If either $p$ is proper, or $A$ is a field, this morphism is an equivalence.

\end{enumerate}

\begin{rmk}\label{realsheaf} The $\ell$-adic realization functor of Definition \ref{d1} can be easily ``sheafified'', as follows. For $X$ a separated, excellent scheme, and $T$ a sheaf of dg-categories on $X_{\textrm{Zar}}$ (i.e. $T \in \dgcat_X := \textrm{lim}\,\dgcat_{\mathcal{U}^{\bullet}}$ where $\mathcal{U}^{\bullet}$ is the Cech-nerve of any Zariski open cover $\mathcal{U}$ of $X$), by taking the limit along the Cech-nerve of a Zariski open cover of $X$, we define $\mathsf{M}_{X}^{-}: \dgcat_X \to  \mathbf{Mod}_{\SH_X}(\BU_X)$ and $\mathsf{R}_{\ell, X} : \mathbf{Mod}_{\SH_X}(\BU_X) \longrightarrow \mathbf{Mod}_{Sh_{\Ql}(X)}(\mathbb{Q}_{\ell, X}(\beta))$. The composite 
$$r_{\ell,X}:= \mathsf{R}_{\ell, X} \circ \mathsf{M}_{X}^{-} : \dgcat_X \longrightarrow \mathbf{Mod}_{Sh_{\Ql}(X)}(\mathbb{Q}_{\ell, X}(\beta))$$ is  called the \emph{$\ell$-adic realization functor for dg-categories over $X$}. Obviously we have $\rl = r_{\ell, S}$, for $\rl$ as in Definition \ref{d1}. Moreover, for $X$ as above, any separated morphism of finite type  $p : X \longrightarrow S$ induces a commutative diagram\footnote{Note that the functor $p_*: \dgcat_X \to \dgcat_S$ consists simply in viewing a dg-category over $X$ as an $A-$dg-category via the morphism $p$. } 
$$\xymatrix{\dgcat_X \ar[d]_{p_*} \ar[r]^-{r_{\ell,X}} & \mathbf{Mod}_{Sh_{\Ql}(X)}(\mathbb{Q}_{\ell, X}(\beta)) \ar[d]^-{p_*} \\ \dgcat_S \ar[r]_-{r_{\ell}} & \mathbf{Mod}_{Sh_{\Ql}(S)}(\mathbb{Q}_{\ell}(\beta))}$$
\end{rmk}

\subsection{The $\ell$-adic Chern character}\label{chernchar}

As explained in \cite{brtv}, there is a symmetric monoidal $\s$-category $\SH_A^{nc}$ of non-commutative motives over $A$.
As an $\s$-category it is the full sub-$\s$-category of $\s$-functors
of (co)presheaves of spectra
$$\dgcat_A^{ft} \longrightarrow \mathbf{Sp},$$
satisfying Nisnevich descent and $\mathbb{A}^1$-homotopy invariance. 

We consider 
$\Gamma : Sh_{\Ql}(S) \longrightarrow \dg_{\Ql},$
the global section $\s$-functor, taking 
an $\ell$-adic complex on $S_{\acute{e}t}$ to its hyper-cohomology. Composing this 
with the Dold-Kan construction $\Map_{\dg_{\Ql}}(\Ql,-) : \dg_{\Ql} \longrightarrow \mathbf{Sp}$  
we obtain an $\s$-functor
$$|-| : Sh_{\Ql}(S) \longrightarrow \mathbf{Sp},$$
which computes hyper-cohomology of $S_{\acute{e}t}$ with $\ell$-adic coefficients, i.e. 
for any $E \in Sh_{\Ql}(S)$, we have natural isomorphisms
$$H^i(S_{\acute{e}t},E) \simeq \pi_{-i}(|E|)\,, i\in \mathbb{Z}.$$
By what we have seen in our last paragraph, the composite functor 
$T \mapsto |\rl(T)|$ provides a
(co)presheaf of spectra
$$\dgcat_A^{ft} \longrightarrow \mathbf{Sp},$$
satisfying Nisnevich descent and $\mathbb{A}^1$-homotopy invariance.
It thus defines an object $|\rl| \in \SH_A^{nc}$. The fact that $\rl$
is lax symmetric monoidal implies moreover that $|\rl|$ is endowed with 
a natural structure of a $\E_\infty$-ring object in $\SH_A^{nc}$. 

Each $T \in \dgcat_A^{ft}$ defines a corepresentable
object $h^T \in \SH_A^{nc}$, characterized by the ($\s$-)functorial equivalence
$$\Map_{\SH_A^{nc}}(h^T,F) \simeq F(T),$$
for any $F \in \SH_A^{nc}$. The existence of $h^T$ is a formal statement, however
the main theorem of \cite{ro} implies that we have a natural equivalence of spectra
$$\Map_{\SH_A^{nc}}(h^T,h^A) \simeq \mathsf{HK}(T),$$
where $\mathsf{HK}(T)$ stands for non-connective homotopy invariant algebraic $K$-theory of the dg-category $T$.
In other words, $T \mapsto \mathsf{HK}(T)$ defines an object in $\SH_A^{nc}$ which 
is isomorphic to $h^B$.
By Yoneda lemma, we thus obtain an equivalence of spaces
$$\Map^{lax-\otimes}(\mathsf{HK},|r_\ell|) \simeq \Map_{\E_\infty-\mathbf{Sp}}(\mathbb{S},|\rl(A)|)\simeq *.$$
In other words, there exists a unique (up to a contractible space of choices) lax symmetric monoidal natural transformation 
$$\mathsf{HK} \longrightarrow |r_\ell|,$$
between lax monoidal $\s$-functors from $\dgcat^{ft}_A$ to $\mathbf{Sp}$. We 
extend this to all dg-categories over $A$, as usual, by passing to Ind-completion
$\dgcat_A \simeq \mathbf{Ind}(\dgcat_A^{ft})$.

\begin{df}\label{d3}
The natural transformation defined above is called the \emph{$\ell$-adic Chern character}.
It is denoted by 
$$Ch_{\ell} : \mathsf{HK}(-) \longrightarrow |\rl(-)|.$$ 
\end{df}

\begin{rmk}Definition \ref{d3} contains a built-in, formal Grothendieck-Riemann-Roch formula. Indeed, 
for any $B$-linear dg-functor $f : T \longrightarrow T'$, the square of spectra
$$\xymatrix{
\mathsf{HK}(T) \ar[r]^-{f_!} \ar[d]_-{Ch_{\ell, T}} & \mathsf{HK}(T') \ar[d]^-{Ch_{\ell, T'}} \\
|\rl(T)| \ar[r]_-{f_!} & |\rl(T')|}$$
commutes up to a natural equivalence. \end{rmk}

\subsection{Trace formula for dg-categories}

Let $\mathcal{C}^{\otimes}$ be a symmetric monoidal $\infty$-category (\cite{tove}, 
\cite[Definition 2.0.0.7]{luha}). 

\begin{hyp}\label{good}
The underlying $\infty$-category $\mathcal{C}$ has small sifted colimits, and the tensor 
product preserves small colimits in each variable. 
\end{hyp}

\begin{df}
Let $\mathcal{C}^{\otimes}$ be a symmetric monoidal $\infty$-category satisfying Hypothesis \ref{good}. We denote by  $\underline{\mathbf{Alg}}(\mathcal{C})$ the $(\infty, 
2)$-category of algebras in $\mathcal{C}^{\otimes}$ denoted by $\mathbf{Alg}_{(1)}
(\mathcal{C}^{\otimes})$ in \emph{\cite[Definition 4.1.11]{lucob}}.
\end{df}

\medskip

Informally, one can describe  $\underline{\mathbf{Alg}}(\mathcal{C})$ as the $(\infty, 
2)$-category with:
\begin{itemize}
\item objects: associative unital monoids (=: $E_1$-algebras) in $\mathcal{C}$.
\item $Map_{\underline{\mathbf{Alg}}(\mathcal{C})}(B,B'):=\mathbf{Bimod}_{B', B}
(\mathcal{C}^{\otimes})$, the $\infty$-category of $(B',B)$-bimodules.  
\item The composition of $1$-morphisms (i.e. of bimodules) is given by tensor product.
\item The composition of $2$-morphisms (i.e. of morphisms between bimodules) is the 
usual composition.
\end{itemize}

\begin{df}\label{rightdual} 
Let $B$ be an algebra in $\mathcal{C}$ and $X$ a left $B$-module. Let us identify $X$ with a 1-
morphism $\underline{X}: 1_{\mathcal{C}} \to B$ in $\underline{\mathbf{Alg}}
(\mathcal{C})$. A \emph{right $B$-dual} of $X$ is defined as a right adjoint $\underline{Y}: B \to 
1_{\mathcal{C}}$ to $\underline{X}$.
\end{df}

Unraveling the definition, we get that a right dual of $X$ is a left $B^{\otimes
\textrm{-op}}$-module $Y$, the unit $u$ of adjunction (or \emph{coevaluation}) is a map 
$coev: 1_{\mathcal{C}} \to Y\otimes_B X $ in $\mathcal{C}$, the counit $v$ of adjunction (or 
\emph{evaluation}) is a map $ev: X\otimes Y \to B$  of $(B,B)$-bimodules; 
$u$ and $v$ satisfy usual compatibilities.\\
Note that, if a right $B$-dual of $X$ exists then it is ``unique'' (i.e. unique up to a 
contractible space of choices).\\

If the right $B$-module $Y$  is the right $B$-dual of the left $B$-module $X$, then we 
can define the trace of any map $f: X\to X$ of left $B$-modules, as follows. \\
\medskip
Recall that we have a coevaluation map in $\mathcal{C}$ $coev: 1_{\mathcal{C}} \to Y
\otimes_B X $ in $\mathcal{C}$ and an evaluation map of $(B,B)$-bimodules $ev: X\otimes 
Y \to B$.\\
Consider the graph $\Gamma_f$ defined as the composite $$\xymatrix{1_{\mathcal{C}} 
\ar[r]^-{coev} & Y\otimes_B X \ar[rr]^-{id\otimes f} & & Y\otimes_B X.}$$
We now elaborate on the evaluation map. Observe that
\begin{itemize}
\item $B \in \mathcal{C}$ has a left $B\otimes B^{\otimes\textrm{-op}}$-module 
structure that we will denote by $B^{L}$.
\item $B \in \mathcal{C}$ has a right $B\otimes B^{\otimes\textrm{-op}}$-module 
structure (i.e. a left $B^{\otimes\textrm{-op}} \otimes B$-module), that we will denote 
by $B^{R}$.
\item $ev: X \otimes Y \to B^{L}$ is a map of left $B\otimes B^{\otimes\textrm{-op}}$-
modules.
\item the composite $\xymatrix{ev': Y\otimes X \ar[r]^-{\sigma} & X\otimes Y \ar[r]^-
{ev} & B^R}$ is a map of left $(B^{\otimes\textrm{-op}}\otimes B)$-modules
\end{itemize}

Apply $(-) \otimes_{B \otimes B^{\otimes\textrm{-op}}}B^L$ to the composite $$
\xymatrix{ev': Y\otimes X \ar[r]^-{\sigma} & X\otimes Y \ar[r]^-{ev} & B^R}$$\\
 to get $$ ev_{HH}: \xymatrix{(Y\otimes X) \otimes_{B \otimes B^{\otimes\textrm{-op}}}B^L  
 \ar[r] & B^R\otimes_{B \otimes B^{\otimes\textrm{-op}}}B^L =: \mathsf{HH}_{\mathcal{C}}
 (B)}.$$
Now observe that $(Y\otimes X) \otimes_{B \otimes B^{\otimes\textrm{-op}}}B^L  \simeq Y
\otimes_B X$ in $\mathcal{C}$. \\Note that, by definition, $\mathsf{HH}_{\mathcal{C}}
 (B)$ is (only) an object in $\mathcal{C}$, called the \emph{Hochschild homology object} of $B$.

\begin{df}\label{nctrace} The \emph{non-commutative trace of} $f: X \to X$ \emph{over $B
$} is defined as the composite
$$Tr_B(f) : \xymatrix{1_{\mathcal{C}} \ar[r]^-{\Gamma_f} & Y\otimes_B X \simeq (Y\otimes 
X) \otimes_{B \otimes B^{\otimes\textrm{-op}}}B^L  \ar[r]^-{ev_{HH}} & B^R 
\otimes_{B \otimes B^{\otimes\textrm{-op}}}B^L =: \mathsf{HH}_{\mathcal{C}}(B)}.$$
$Tr_B(f)$ is a morphism in $\mathcal{C}$.
\end{df}

\begin{rmk}\label{tr-commvsnoncomm} Let $B \in \mathbf{CAlg}(\mathcal{C^{\otimes}})$, 
and let us still denote by $B$ its image via the canonical map $\mathbf{CAlg}
(\mathcal{C}^{\otimes}) \to \Alg_{E_1}(\mathcal{C}^{\otimes})$. In this case, $
\mathbf{Mod}_B(\mathcal{C}^{\otimes})$ is a symmetric monoidal $\infty$-category, and if 
$X\in \mathbf{Mod}_B(\mathcal{C}^{\otimes})$ is a dualizable object (in the usual 
sense), then its (left and right) dual in $\mathbf{Mod}_B(\mathcal{C}^{\otimes})$ is 
also a right-dual of $X$ according to Definition \ref{rightdual}. In this case, any $f: X \to X$ in $\mathbf{Mod}_B(\mathcal{C}^{\otimes})$,  has therefore \emph{two possible traces}, a non-commutative one (as in Definition \ref{nctrace}) $$Tr_B(f): 1_{\mathcal{C}} \longrightarrow \mathsf{HH}_{\mathcal{C}}(B)$$ which is a morphism in $\mathcal{C}$, and a more standard, commutative one $$Tr^{\textrm{c}}_B(f): B \longrightarrow B$$ which is a morphism on $\mathbf{Mod}_B(\mathcal{C}^{\otimes})$. The two traces are related by the following commutative diagram $$\xymatrix{1_{\mathcal{C}} \ar[r]^-{u_B} \ar[d]_-{Tr_B(f)} & B \ar[d]^-{Tr^{\textrm{c}}_B(f)} \\ \mathsf{HH}_{\mathcal{C}}(B) \ar[r]_-{a} & B }$$ where $a: \mathsf{HH}_{\mathcal{C}}(B) \to B$ is the canonical augmentation (which exists since $B$ is commutative), and $u_B: 1_{\mathcal{C}} \to B$ is the unit map of the algebra $B$ in $\mathcal{C}$.  
\end{rmk}

\bigskip

\noindent \textbf{The case of dg-categories.} Let us specialize the previous discussion 
to the case $\mathcal{C}^{\otimes}= \mathbf{dgCat}_A$.

\noindent Let $B$ be a monoidal dg-category, i.e. an associative and unital monoid in the 
symmetric monoidal $\infty$-category $\mathbf{dgCat}_A$.

\begin{prop}\label{prop-alwaysduals}
For any $B$-module $T$ which is cotensored in the sense of Definition \ref{dcotens}, the big
$B$-module $\widehat{T} \in \dgCAT_B$ has a right dual in the symmetric
monoidal $\infty$-category $\dgCAT_A$ whose underlying 
big dg-category is $\widehat{T^o}$. 
\end{prop}

\noindent \textbf{Proof.} This is very similar to the argument used in 
\cite[Prop. 2.5 (1)]{to2}. We consider $\widehat{T^o}$, and we define
evaluation and coevaluation maps as follows. 

The big morphism $h$ introduced right before Definition \ref{dcotens}
$$h : \widehat{T^o} \widehat{\otimes}_A \widehat{T} \longrightarrow \widehat{B},$$
whose domain is naturally a $\widehat{B}$-bimodule, can be canonically lifted to a morphism of bimodules. We choose this as our evaluation morphism.

The coevaluation is then obtained by duality. We start by 
the diagonal bimodule
$$T : (T^o \otimes_A T)^o \longrightarrow L(A)=\widehat{A}.$$
sending $(x,y)$ to $T(y,x)$. This morphism naturally descends
to $(T^o \otimes_B T)^o$, providing a dg-functor
$$(T^o \otimes_B T)^o\longrightarrow \widehat{A}.$$
Note that $T^o$ is naturally a right $B$-module since $T$ is assumed to be
cotensored (and note that, otherwise, $T^o \otimes_B T$ would not make sense).
This dg-functor is an object in $\widehat{T^o\otimes_B T} \simeq
\widehat{T^o} \widehat{\otimes}_A \widehat{T}$, and thus defines
a coevaluation morphism
$$\widehat{A} \longrightarrow \widehat{T^o} \widehat{\otimes}_A \widehat{T}.$$
These two evaluation and coevaluation morphisms satisfy the 
required triangular identities, and thus make $\widehat{T^o}$ a right
dual to $\widehat{T}$. \hfill $\Box$\\

According to the previous Proposition, any cotensored
$B$-module $T$ has a big right dual, so it comes equipped with 
big evaluation and coevaluation maps.

\begin{df}
For a monoidal dg-category $B$, and a $B$-module
$T \in \mathbf{dgCat}_B$, we say that $T$ is \emph{saturated over} $B$ if 

\begin{enumerate}

\item $T$ is cotensored (over $B$), and

\item the evaluation and coevaluation morphisms are small (i.e. are maps in $\dgcat_A$).

\end{enumerate}

\end{df}

In particular, if $T$ saturated over $B$, and $f:T \to T$ is a morphism in $
\mathbf{dgCat}_\B$,
then the trace $$Tr_{B}(f:T \to T): A \to \mathsf{HH}(B/A)=B^R \otimes_{B^{\otimes
\textrm{-op}}\otimes_A B}B^L $$ is also small, i.e. it is a morphism \emph{inside}  $
\mathbf{dgCat}_A$.\\
\medskip

Since our $\ell$-adic realization functor $\rl$ is only lax-monoidal, in order to establish our trace formula for an endomorphism $f: T\to T$ of a saturated $B$-module, we need to restrict to those saturated $B$-modules on which $\rl$ is in fact symmetric monoidal.
 
\begin{df}\label{ammisibilta}
A saturated $T \in \mathbf{dgCat}_{B}$ is called \emph{$\ell^{\otimes}$-admissible}  if the canonical 
map $$\rl(T^{op})\otimes_{\rl(B)} \rl(T) \to \rl(T^{op}\otimes_{B} T)$$ is an 
equivalence in $Sh_{\ql}(S)$. 
\end{df}

The trace $Tr_B(f)$ is a map $A \to \mathsf{HH}(B/A)$, hence it induces a map in $\mathbf{Sp}$ $$K(Tr_B (f)): K(A) \to K(\mathsf{HH}(B/A))$$  which is actually a map of $K(A)$-modules (in spectra), since $K$ is lax-monoidal. Hence it corresponds to an element denoted as $$ [\mathsf{HH}(T/B,f)] \equiv tr_B (f) \in K_0(\mathsf{HH}(B/A)).$$ Therefore, its image by the $\ell$-adic Chern character $Ch_{\ell, 0}:= \pi_0(Ch_{\ell})$ $$Ch_{\ell, 0}:  K_0(\mathsf{HH}(B/A)) \to Hom_{D(\rl(A))}(\rl(A), \rl(\mathsf{HH}(B/A))) \simeq H^0(S_{\textrm{\'et}}, \rl(\mathsf{HH}(B/A))),$$ is an element $Ch_{\ell, 0}([\mathsf{HH}(T/B;f)]) \in  H^0(S_{\textrm{\'et}}, \rl(\mathsf{HH}(B/A)))$.

\

On the other hand, the trace $Tr_{\rl(B)}(\rl(f))$ of $\rl(f)$ over $\rl(B)$ is, by definition, a morphism $$\rl(A) \simeq \mathbb{Q}_{\ell}(\beta) 
\to  \mathsf{HH}(\rl(B)/ \rl(A))$$ in $\mathbf{Mod}_{\rl(A)}(Sh_{\mathbb{Q}_{\ell}}(S))$. \\ We may further compose this with the canonical map $$\mathsf{HH}(\rl(B)/ \rl(A)) \to \rl(\mathsf{HH}(B/ A))$$ (given by lax-monoidality of $\rl(-)$), to get a map  $$\rl(A) \to \rl(\mathsf{HH}(B/ A))$$ in $\mathbf{Mod}_{\rl(A)}(Sh_{\mathbb{Q}_{\ell}}(S))$. This is the same thing as an element denoted as 
$$tr_{\rl(B)}(\rl(f)) \in \pi_0 (|\rl(\mathsf{HH}(B/ A))|) \simeq H^0(S_{\textrm{\'et}}, \rl(\mathsf{HH}(B/ A)) ).$$

\begin{thm}\label{ttrace}
Let $B$ a monoidal dg-category over $A$, $T\in \mathbf{dgCat}_B$ a saturated and $
\ell^{\otimes}$-admissible $B$-module, and $f:T\to T$ map in $\mathbf{dgCat}_B$. Then,
we have an equality 
$$Ch_{\ell, 0}([\mathsf{HH}(T/B,f)]) =tr_{\rl(B)}(\rl(f)) $$ in $H^0(S_{\textrm{\'et}}, \rl(\mathsf{HH}
(B/ A))$. 
\end{thm}

\noindent \textbf{Proof.} This is a formal consequence 
of uniqueness of right duals, and of the resulting fact that traces are preserved
by symmetric or lax symmetric monoidal $\s$-functors under our admissibility condition.  The key statement 
is the following lemma, left as an exercise to the reader, and 
applied to the case where $F$ is our $\ell$-adic realization functor.

\begin{lem}\label{ltrace}
Let $$F : \mathcal{C} \longrightarrow \D$$ be a lax symmetric monoidal $\s$-functor between  
presentable symmetric monoidal $\s$-categories.
Let $B$ be a monoid in $\mathcal{C}$, $M$ a left $B$-module,  
and $f : M \longrightarrow M$ a morphism of $B$-modules. 
We assume that $M$ has a right dual $M^o$, and that the natural
morphism
$$F(M^o) \otimes_{F(B)}F(M) \longrightarrow F(M^o\otimes_B M)$$
is an equivalence. Then, $F(M)$ has a right dual, and we have 
$$F(Tr(f)) = i(Tr(F(f)))$$
as elements in $\pi_0(Hom_{\D}(\mathbf{1},F(\mathsf{HH}(B)))$, where $i$
is induced by the natural morphism $\mathsf{HH}(F(B)) \longrightarrow F(\mathsf{HH}(B))$.
\end{lem}

\hfill $\Box$ \\

\section{Invariant vanishing cycles}

This section gathers general results about inertia-invariant vanishing cycles 
(\textit{I-vanishing cycles}, for short), their relations with dg-categories
of singularities, and their behaviour under products. These results are partially taken
from \cite{brtv}, and the only original result is Proposition \ref{kunn} that can be seen as 
a version of Thom-Sebastiani formula in the mixed-characteristic setting. \\

All along this section, $A$ will be a strictly henselian excellent dvr with fraction field $K=\mathsf{Frac}(A)$, and
perfect (hence algebraically closed) residue field $k$. We let $S=Spec\, A$. All schemes
over $S$ are assumed to be separated and of finite type over $S$.
We denoted by $i : s:=Spec\, k \longrightarrow S$ the closed point of $S$, and
$j : \eta := Spec\, K \longrightarrow S$ its generic point. 
For an $S$-scheme $X$, we denote by $X_s:=X\times_S s$ its special fiber, and 
$X_\eta=X\times_S \eta$ its generic fiber. Accordingly, we write $X_{\bar{\eta}}:=X\times_S Spec\, K^{sp}$
for the \emph{geometric} generic fiber. \\

\subsection{Trivializing the Tate twist}\label{trivtate}We let $\ell$ be a prime invertible in $k$, and we denote by $p$ the characteristic exponent of $k$. 
As $k$ is algebraically closed, we may, and will, choose
once for all a group isomorphism 
$$\mu_{\infty}(k)\simeq \mu_{\infty}(K) \simeq (\mathbb{Q}/\mathbb{Z})[p^{-1}]$$
between the group of roots of unity in $k$ and the prime-to-$p$ part of
$\mathbb{Q}/\mathbb{Z}$. Equivalently, we have chosen a given group isomorphism 
$$\lim_{(n,p)=1}\mu_n(k) \simeq \hat{\mathbb{Z}}',$$
where $\hat{\mathbb{Z}}':=\lim_{(n,p)=1} \mathbb{Z}/n$. In particular, we have selected a 
topological generator of $\lim_{(n,p)=1}\mu_n(k)$, corresponding to the image
of $1\in \mathbb{Z}$ inside $\hat{\mathbb{Z}}'$. 
The choice of the isomorphism above also provides a specific induced 
isomorphism $\Ql(1) \simeq \Ql$ of $\Ql$-sheaves on $S$, where $(1)$ denotes, as usual, the Tate twist. By taking tensor powers of this isomorphism, we get 
various induced isomorphisms $\Ql(i) \simeq \Ql$ for all $i\in \mathbb{Z}$. 

We remind that the absolute Galois group $\mathrm{I}$ of $K$ (which coincides with the \emph{inertia group} in our case) sits in an extension of pro-finite
groups\footnote{Note that the tame inertia quotient $\mathrm{I}_t$ is canonically isomorphic to $\hat{\mathbb{Z}}'(1)$, and it becomes isomorphic to $\hat{\mathbb{Z}}'$ through our choice.  }
$$\xymatrix{
1 \ar[r] & P \ar[r] & \mathrm{I} \ar[r] & \mathrm{I}_t \simeq \hat{\mathbb{Z}}' \ar[r] & 1\, ,}$$
where $P$ is a pro-p-group (the \emph{wild inertia} subgroup). For any continuous finite dimensional $\Ql$-representation
$V$ of $\mathrm{I}$, the group $P$ acts by a finite quotient $G_V$ on $V$. Moreover, 
the Galois cohomology of $V$ can then be explicitly identified with the two-terms complex
$$\xymatrix{
V^G \ar[r]^-{1-T} & V^G}$$
where $T$ is the action of the chosen topological generator of $\mathrm{I}_t$. This
easily implies that for any $\Ql$-representation $V$ of $\mathrm{I}$, the natural pairing
on Galois cohomology
$$H^i(\mathrm{I},V) \otimes H^{1-i}(\mathrm{I},V^{\vee}) \longrightarrow H^{1}(\mathrm{I},\Ql)\simeq \Ql$$
is non-degenerate. In other words, if we denote by $V^\mathrm{I}$ the complex
of cohomology of $\mathrm{I}$ with coefficients in $V$, we have a natural quasi-isomorphism
$(V^\mathrm{I})^{\vee} \simeq (V^{\vee})^\mathrm{I}[1].$

\subsection{Reminders on actions of the inertia group}

Let $X \longrightarrow S$ be an $S$-scheme (separated and of finite type, according to our conventions). 
We recall from \cite[Exp. XIII, 1.2]{sga7II} that 
we can associate to $X$ a \emph{vanishing topos} $(X/S)_{et}^{\nu}$ which is defined as
(a 2-)fiber product of toposes
$$(X/S)_{et}^{\nu}:=(X_s)^{\sim}_{et} \times_{s_{et}^{\sim}}\eta_{et}^{\sim}.$$
Since $S$ is strictly henselian, $s_{et}^{\sim}$ is  in fact the punctual topos, and the
fiber product above is in fact a product of topos. The topos 
$\eta_{et}^{\sim}$ is equivalent to the topos of sets with continuous 
action of $\mathrm{I}=\Gal(K^{sp}/k)$, where $K^{sp}$ denotes a seprable closure of $K$. 
Morally, $(X/S)_{et}^{\nu}$ is the topos
of \'etale sheaves on $X_s$, endowed with a continuous action of $\mathrm{I}$ (see \cite[Exp. XIII, 1.2.4]{sga7II}).

As explained in \cite{brtv} we have an $\ell$-adic $\s$-category 
$\D((X/S)_{et}^{\nu},\mathbb{Z}_\ell)$. Morally speaking, objects of this $\s$-category
consist of the data of an object $E \in \D(\bar{X}_s,\mathbb{Z}_\ell)= Sh_{\mathbb{Z}_{\ell}}(\bar{X}_s)$ together 
with a continuous action of $\mathrm{I}$. We say that such 
an object is \emph{constructible} if $E$ is a constructible object in 
$\D(\bar{X}_s,\mathbb{Z}_\ell)=Sh_{\mathbb{Z}_{\ell}}(\bar{X}_s)$, and we denote by $\D_c((X/S)_{et}^{\nu},\mathbb{Z}_\ell)$
the full sub-$\s$-category of constructible objects. 

\begin{df}
The \emph{$\s$-category of ind-constructible 
$\mathrm{I}$-equivariant $\ell$-adic complexes on $X_s$} is defined by 
$$\D^{\mathrm{I}}_{ic}(X_s,\Ql):=\Ind(\D_{c}((X/S)_{et}^{\nu},\mathbb{Z}_\ell) \otimes_{\mathbb{Z}_\ell}
\Ql).$$
The full sub-$\s$-category of constructible objects is 
$\D^{\mathrm{I}}_{c}(X_s,\Ql):=\D_{c}((X/S)_{et}^{\nu},\mathbb{Z}_\ell) \otimes_{\mathbb{Z}_\ell}
\Ql.$
\end{df}
 
Note that since 
we have chosen trivialisations of the Tate twists, $\Ql(\beta)$ is identified with
$\Ql[\beta,\beta^{-1}]$ where $\beta$ is  a free variable in degree $2$. This is a graded algebra object in 
$\D^{\mathrm{I}}_{ic}(X_s,\Ql)$, and we define 
the $\s$-category $\D^{\mathrm{I}}_{ic}(X_s,\Ql(\beta))$ as the $\s$-category of $\Ql(\beta)$-modules in $\D^{\mathrm{I}}_{ic}(X_s,\Ql)$, or
equivalently, the $2$-periodic $\s$-category of ind-constructible $\Ql$-adic
complexes on $(X/S)^{\nu}_{et}$.

\begin{df}
An object $E \in \D^{\mathrm{I}}_{ic}(X_s,\Ql(\beta))$ is \emph{constructible} if
it belongs to the thick triangulated sub-$\s$-category generated by objects
of the form $E_0(\beta)=E_0\otimes_{\Ql}\Ql(\beta)$ for $E_0$ a constructible
object in $\D_{ic}^\mathrm{I}(X_s,\Ql)$. 

The full sub-$\s$-category of $\D^{\mathrm{I}}_{ic}(X_s,\Ql(\beta))$ consisting of constructible
objects is denoted $\D^{\mathrm{I}}_{c}(X_s,\Ql(\beta))$. Similarly, for any $S$-scheme $X$, 
we define $\D_c(X,\Ql(\beta))$ as the full sub-$\s$-category of 
objects in $\D_{ic}(X,\Ql(\beta))$ generated by $E_0(\beta)$ for 
$E_0$ constructible.
\end{df}

Note that strictly speaking an object of $\D^{\mathrm{I}}_{c}(X_s,\Ql(\beta))$ is not
constructible in the usual sense, as its underlying object 
in $\D_{ic}(X_s,\Ql)$ is $2$-periodic. \\

The topos $(X/S)_{et}^{\nu}$ comes with a natural projection
$(X/S)^{\nu}_{et} \longrightarrow (X_s)^{\sim}_{et}$ whose direct image is an $\s$-functor denoted by
$$(-)^\mathrm{I} : \D^{\mathrm{I}}_{ic}(X_s,\Ql) \longrightarrow \D_{ic}(X_s,\Ql)$$
called the \emph{$\mathrm{I}$-invariants} functor. This $\s$-functor preserves 
constructibility. The $\s$-categories $\D^{\mathrm{I}}_{ic}(X_s,\Ql)$ and $\D_{ic}(X_s,\Ql)$ 
carries natural symmetric monoidal structures and the $\s$-functor 
$(-)^\mathrm{I}$ comes equipped with a natural lax symmetric monoidal structure (being 
induced by the direct image of a morphism of toposes). Moreover, 
$(-)^\mathrm{I}$ is the right adjoint of the symmetric monoidal $\s$-functor
$U: \D_{ic}(X_s,\Ql) \longrightarrow \D^{\mathrm{I}}_{ic}(X_s,\Ql)$
endowing objects in $\D_{ic}(X_s,\Ql)$ with the trivial
action of $\mathrm{I}$. This gives $\D^{\mathrm{I}}_{ic}(X_s,\Ql)$
the structure of a $\D_{ic}(X_s,\Ql)$-module via $\D_{ic}(X_s,\Ql) \times \D^{\mathrm{I}}_{ic}(X_s,\Ql) \to \D^{\mathrm{I}}_{ic}(X_s,\Ql) : (E,F) \mapsto U(E)\otimes F $. This bi-functor distributes over colimits, thus by the adjoint theorem, we get an enrichment of $\D^{\mathrm{I}}_{ic}(X_s,\Ql)$ over $\D_{ic}(X_s,\Ql)$. Note that $\D^{\mathrm{I}}_{ic}(X_s,\Ql)$ is also enriched over itself.\\

It is important to notice that the $\mathrm{I}$-invariants functor commutes with base change
in the following sense. Let $f : \mathrm{Spec}\, k \longrightarrow X_s$ be  a 
geometric point. The morphism $f$ defines a geometric point of $(X_s)_{et}^{\sim}$ and thus induces 
a geometric morphism of toposes
$$\eta_{et}^{\sim} \longrightarrow (X/S)^{\nu}_{et}.$$
We thus have an inverse image functor
$$f^* : \D_{c}^{\mathrm{I}}(X_s,\Ql) \longrightarrow 
\D_c(\eta,\Ql)=D_{c}(\mathrm{I}, \Ql).$$
where $\D_c(\mathrm{I}, \Ql)$ is the $\s$-category of finite dimensional complexes of $\ell$-adic
representations of $\mathrm{I}$.
As usual, the square of $\s$-functors 
$$\xymatrix{
\D_c^\mathrm{I}(X_s,\Ql) \ar[r]^-{(-)^\mathrm{I}} \ar[d]_-{f^*} & \D_c(X_s,\Ql)=Sh^{\mathsf{ct}}_{\Ql}(X_s) \ar[d]^-{f^*} \\
\D_c(\mathrm{I},\Ql) \ar[r]_-{(-)^\mathrm{I}} & \D_c(\Ql)}$$
comes equipped with a natural transformation
$$f^*((-)^\mathrm{I}) \Rightarrow (f^*(-))^\mathrm{I}.$$
It can be checked that this natural transformation is always an 
equivalence. In particular, for any geometric point $x$ in $X_s$, 
we have a natural equivalence
of $\ell$-adic complexes
$(E)^{\mathrm{I}}_x \simeq (E_x)^\mathrm{I}$, 
for any $E \in \D_c^\mathrm{I}(X_s,\Ql)$. \\

The dualizing complex $\omega$ of the scheme $X_s$ can be used in order
to obtain a dualizing object in $\D_{c}^\mathrm{I}(X_s,\Ql)$ as follows. 
We consider $\omega$ as an object in $\D_c^\mathrm{I}(X_s,\Ql)$ endowed
with the trivial $\mathrm{I}$-action. We then have an 
equivalence of $\s$-categories
$$\mathbb{D}_I : \D_c^\mathrm{I}(X_s,\Ql) \longrightarrow \D_c^\mathrm{I}(X_s,\Ql)^{op}$$
sending $E$ to $\mathbb{R}\underline{Hom}(E,\omega)$, 
where $\mathbb{R}\underline{Hom}$ denotes the natural enrichment of 
$\D_c^{\mathrm{I}}(X_s,\Ql)$ over itself. 
We obviously have a canonical biduality 
equivalence $\mathbb{D}_{\mathrm{I}}^2 \simeq id$. The duality functor $\mathbb{D}_I$ is compatible with
the usual Grothendieck duality functor $\mathbb{D}$ for the scheme $X_s$ up to a shift, as explained by the following 
lemma.

\begin{lem}\label{I-duality}
For any object $E\in \D_{c}^\mathrm{I}(X_s,\Ql)$, there is a functorial equivalence in 
$\D_{c}(X_s,\Ql)$
$$d : \mathbb{D}(E^\mathrm{I})[-1]\simeq (\mathbb{D}_I(E))^\mathrm{I}.$$
\end{lem}

\noindent \textbf{Proof.} Taking $\mathrm{I}$-invariants is a lax monoidal $\s$-functor, so we have a natural
map
$E^\mathrm{I} \otimes \mathbb{D}_I(E)^\mathrm{I} \longrightarrow (E\otimes \mathbb{D}_I(E))^\mathrm{I},$
that can be composed with the evaluation morphism
$E\otimes \mathbb{D}_I(E) \longrightarrow \omega$ to obtain
$E^\mathrm{I} \otimes \mathbb{D}_I(E)^\mathrm{I} \longrightarrow \omega^\mathrm{I}.$
As the action of $\mathrm{I}$ on $\omega$ is trivial, we have a canonical equivalence 
$\omega^\mathrm{I}\simeq \omega \otimes \Ql^\mathrm{I} \simeq \omega \oplus \omega[-1]$. 
By projection on the second factor we get a pairing
$E^\mathrm{I}\otimes \mathbb{D}_I(E)^\mathrm{I} \longrightarrow \omega[-1]$,
and thus a map
$$\mathbb{D}_I(E)^\mathrm{I} \longrightarrow  \mathbb{D}(E^\mathrm{I})[-1].$$
We claim that the above morphism is an equivalence in $\D_c(X_s,\Ql)$. 
For this it is enough to check that the above morphism is a stalkwise equivalence. 
Now, the stalk of the above morphism at a geometric point $x$ of $X_s$ can be written as 
$$(E(x)^{\vee})^\mathrm{I} \longrightarrow (E(x)^\mathrm{I})^{\vee}[-1]$$
where $E(x):=H^*_x(X,E) \in D_c(\Ql)$ is the local cohomology of $E$ at $x$,  
and $(-)^{\vee}$ is now the standard linear duality over $\Ql$. The result now follows from the following well-known 
duality for $\Ql$-representations of $\mathrm{I}$: for any finite dimensional $Ql$-representation 
$V$ of $\mathrm{I}$, the fundamental class in $H^{1}(\mathrm{I},\Ql)\simeq \Ql$, induces 
a canonical isomorphism of Galois cohomologies
$$H^*(\mathrm{I},V^{\vee}) \simeq H^{1-*}(\mathrm{I},V)^{\vee}.$$
\hfill $\Box$ \\

\subsection{Invariant vanishing cycles and dg-categories}

From \cite[Exp. XIII]{sga7II} and \cite[4.1]{brtv}, the vanishing cycles construction provides an $\s$-functor
$$\phi : \D_{c}(X,\Ql) \longrightarrow \D_{c}^{\mathrm{I}}(X_s,\Ql).$$
Applied to the constant sheaf $\Ql$, we get this way an object 
denoted by $\nu_{X/S}$ (or simply $\nu_X$ if $S$ is clear) in 
$\D_{c}^{\mathrm{I}}(X_s,\Ql)$.

\begin{df}
The \emph{$\mathrm{I}$-invariant vanishing cycles} of $X$ relative to $S$ (or \emph{$\mathrm{I}$-vanishing cycles},  
for short) is the object 
$$\nu_X^\mathrm{I}:= (\nu_X)^\mathrm{I} \in \D_{c}(X_s,\Ql).$$
\end{df}

There are several possible descriptions of invariant vanishing cycles. First of all, by its very definition, 
$\nu_X^\mathrm{I}$ is related to the $\mathrm{I}$-invariant nearby cycles $\psi_X^{\mathrm{I}}:= (\psi_X)^{\mathrm{I}}$ by means of an exact triangle
in $D_c(X_s,\Ql)$
\begin{equation}\label{1}\xymatrix{
\Ql^\mathrm{I} \ar[r] & \psi_X^{\mathrm{I}} \ar[r] & \nu_X^\mathrm{I}.}
\end{equation}
Another description,  
in terms of \emph{local cohomology}, is the following. We let $U=X_K$ be the open complement 
of $X_s$ inside $X$, and $j_X : U \hookrightarrow X$ and $i_X : X_s \hookrightarrow X$ the corresponding immersions. 
Then, the $\mathrm{I}$-vanishing cycles enters in an exact triangle in $D_c(X_s,\Ql)$ 
\begin{equation}\label{2}\xymatrix{\nu_X^\mathrm{I} \ar[r] & \Ql \ar[r] & i_X^!(\Ql)[2].}
\end{equation}
Triangle $(2)$ follows from the octahedral axiom applied to the triangles $(1)$ and $$\Ql \longrightarrow \Ql^\mathrm{I} \simeq \Ql \oplus \Ql[-1] \longrightarrow \Ql[-1],$$
taking also into account the triangle $$i_X^! \Ql \longrightarrow \Ql \longrightarrow i_X^*(j_X)_*j_X^*\Ql \simeq \psi_X^\mathrm{I}.$$

We get one more description of $\nu_X^\mathrm{I}$  (or rather, of $\nu_X^\mathrm{I}(\beta):=\nu_X^\mathrm{I}\otimes \Ql(\beta)$) using the \emph{$\ell$-adic realization of the dg-category of singularities} studied 
in \cite{brtv}, at least when $X$ is a regular scheme with smooth generic fiber.
Let $\Sing(X_s)= \Coh (X_s)/\Perf (X_s)$ be the dg-category of singularities of the scheme $X_s$. This
dg-category is naturally linear over the dg-categories $\Perf(X_s)$ and $\Perf(X)$, and thus we can take its $\ell$-adic realization  
 $r_{\ell, X}(\Sing(X_s))$ over $X$ (see Remark \ref{realsheaf}) which is a $\mathbb{Q}_{\ell, X}(\beta)$-module in $\D_{ic}(X,\Ql)$ supported on $X_s$, hence can be identified with a $\mathbb{Q}_{\ell, X_s}(\beta)$-module
in $\D_{ic}(X_s,\Ql)$. When $X$ is a regular scheme and 
$X_K$ is smooth over $K$, we have from \cite{brtv} a canonical 
equivalence in $\D_{c}(X_s,\Ql(\beta))$ \footnote{Strictly speaking, in \cite{brtv} the equivalence (\ref{prelim}) is proved only after push-forward to $S$ but the very same proof shows also the equivalence (\ref{prelim}).}
\begin{equation}\label{prelim}\nu_X^\mathrm{I}(\beta)[1] \simeq r_{\ell, X}(\Sing(X_s))
\end{equation}
where $\nu_X^\mathrm{I}(\beta)$ stands for $\nu_X^\mathrm{I}\otimes \Ql(\beta)$. \\

We conclude this section with another description of $\nu_X^\mathrm{I}(\beta)$, see equivalence (\ref{final}), this time in terms
of \emph{sheaves of singularities}. We need a preliminary result.

\begin{lem}\label{fund} We assume that $S$ is excellent.
Let $p: X \to S$ be a separated morphism of finite type.
We write $r_{\ell, X} : \SH_{X} \to Sh_{\Ql} (X)$, for the $\ell$-adic realization \emph{over $X$}, as in Remark \ref{realsheaf}.Then, we have
\begin{enumerate}
\item $r_{\ell, X}(\Perf (X)) \simeq \mathbb{Q}_{\ell, X} (\beta)$ in $\D_{ic}(X,\Ql)$.
\item $r_{\ell, X}(\Coh (X)) \simeq \omega_X (\beta)$ in $\D_{ic}(X,\Ql)$, where $\omega_X \simeq p^!(\Ql)$ is the $\ell$-adic dualizing complex of $X$\footnote{Note that since $S$ is excellent, $X$ is excellent so that $\omega_X$ exists by a theorem of Gabber (\cite[Exp. XVII, Th. 0.2]{ilo}).}.
\item There exists a canonical map $\eta_X : \mathbb{Q}_{\ell, X}(\beta) \longrightarrow \omega_X(\beta)$ in $\D_c(X,\Ql(\beta))$, called the \emph{2-periodic $\ell$-adic fundamental class} of $X$.
\end{enumerate}
\end{lem}

\noindent \textbf{Proof.} First of all, separated finite type morphisms
of noetherian schemes are compactifiable (by Nagata's theorem), thus we can 
assume that $p: X \to S$ is proper.

(1) follows immediately from \cite[Prop. 3.9 and formula (3.7.13))]{brtv}. In order to prove (2) we first produce a map $\alpha: \rl(\Coh (X)) \to \omega_X (\beta)$. In the notations of \cite[\S 3]{brtv} (note that \cite{brtv}'s notation for $\mathsf{M}_X^{T}$ is $\mathcal{M}_X^{\vee}(T)$), we first construct a map  $\alpha^{\textrm{mot}}: \mathcal{M}^{\vee}_X(\Coh (X)) \to p^!(\BU_S)=: \omega^{\textrm{mot}}_X$ in $\mathrm{SH}(X)$, whose \'etale $\ell$-adic realization will be $\alpha$. Since $p$ is proper, $\alpha^{\textrm{mot}}$ is the same thing, by adjunction, as a map $p_*(\mathcal{M}^{\vee}_X(\Coh (X))) \to \BU_S$ in $\mathrm{SH}(S)$.
Now, $p_*(\mathcal{M}^{\vee}_X(\Coh (X)))$ is just $\mathcal{M}^{\vee}_S(\Coh (X)))$, where $\Coh (X)$ is viewed as a dg-category over $S$, via $p$. If $Y$ is smooth over $S$, we have by \cite[Prop. B.4.1]{pre}, an equivalence of $S$-dg-categories
\begin{equation}\label{eq-cohperf=coh}
\Coh(X)\otimes_S \Coh(Y) \simeq \Coh(X\times_S Y)
\end{equation}
Through this identification, $\mathcal{M}_S^\vee(\Coh(X))\in \mathrm{SH}(S)$ is the $\infty$-functor 
$Y\mapsto \mathrm{KH}(\Coh(X\times_S Y))$, and $\mathrm{KH}(\Coh(X\times_S Y))$  is equivalent to the $\mathrm{G}$-theory spectrum  $\mathrm{G}(X\times_S Y)$ of $X\times_S Y$, by $\mathbb{A}^1$-invariance of G-theory. Since $S$ is regular, $\BU_S \simeq \mathrm{G}_S:= \mathrm{G}(-/S)$ canonically in $\mathrm{SH}(S)$, and we can take the map $\mathcal{M}^{\vee}_S(\Coh (X))) \to \BU_S \simeq \mathrm{G}_S$ to be the push forward $p_*$ on G-theories $\mathrm{G}(X\times_S -) \to G(-/S)$. This gives us a map $\alpha^{\textrm{mot}}: \mathcal{M}^{\vee}_X(\Coh (X)) \to p^!(\BU_S)=: \omega^{\textrm{mot}}_X$. Now observe that by \cite[formula (3.7.13)]{brtv}, the \'etale $\ell$-adic realization of $p^!(\BU_S)$ is canonically equivalent to $p^!(\Ql (\beta)) \simeq \omega_X (\beta)$ (since \'etale $\ell$-adic realization commutes with six operations, \cite[Rmk. 3.23]{brtv}). Therefore, we get our map $\alpha: \rl(\Coh (X)) \to \omega_X (\beta)$. Checking that $\alpha$ is an equivalence is a local statement, i.e. it is enough to show that if $j: V=\mathrm{Spec}, A \hookrightarrow X$ is an open affine subscheme, then $j^*(\alpha)$ is an equivalence. Now, $j^*\rl(\Coh (X)) \simeq r_{\ell, V}(j^*\Coh (X)) \simeq r_{\ell, V}(\Coh (V))$ (where $r_{\ell, V}$ denotes the $\ell$-adic realization over $V$), and $j^*\omega_X \simeq j^!\omega_X \simeq \omega_V$, so $j^*\alpha$ identifies with a map $r_{\ell, V}(\Coh (V))\to \omega_V(\beta)$. Since $V$ is affine and of finite type over $S$, we can choose a closed immersion $i:V\hookrightarrow V'$, with $V'$ affine and smooth  (hence regular) over $S$. Let $h: V'\setminus V \hookrightarrow V'$ be the complementary open immersion. Since $V'$ and $V'\setminus V$ are regular, by Quillen localization and the properties of the nc realization functor $\mathcal{M}^{\vee}$ (see \cite{brtv}), we get a cofiber sequence $$\mathcal{M}^{\vee}_{V'}(\Coh (V)/V') \to \BU_{X'} \to h_*\BU_{V'\setminus V}, $$ where the notation $\Coh (V)/V'$ means that $\Coh (V)$ is viewed as a dg-catgeory over $V'$, via $i$. In other words, $\mathcal{M}^{\vee}_{V'}(\Coh (V)/V')\simeq i_*\mathcal{M}^{\vee}_{V}(\Coh (V))$. If we apply $i^*$ to this cofiber sequence, and compare what we obtain to the application of $i^*$ to the standard localization sequence $$i_*i^! \BU_{V'} \to \BU_{V'}\to h_*h^*\BU_{V'}=h_*\BU_{V'\setminus V},$$ we finally get, after \'etale $\ell$-adic realization, that $\omega_V(\beta) \simeq i^!\Ql (\beta) \simeq r_{\ell, V}(\Coh (V))$. This implies that $j^*\alpha$ is also an equivalence. \\ By (1) and (2), the map in (3) is finally obtained by applying 
$r_{\ell, X}$ to the inclusion $\Perf (X) \to \Coh (X)$.

\hfill $\Box$

\begin{rmk} Note that Lemma \ref{fund} applies, in particular, to give a 2-periodic $\ell$-adic fundamental class map $\eta_U : \Ql(\beta) \longrightarrow \omega_U(\beta)$ for any open subscheme $U\hookrightarrow X$ over $S$, whenever $X$ is proper over $S$.
\end{rmk}

\begin{df}
Let $X/S$ be a separated and finite type $S$-scheme.
The \emph{sheaf of singularities of $X$} is defined to be the cofiber of the 2-periodic 
$\ell$-adic fundamental class morphism $\eta_X$ (Lemma \ref{fund} (3))
$$\omega_X^o := \textrm{Cofib}(\eta_X :\Ql(\beta) \longrightarrow \omega_X(\beta)).$$ 
\end{df}

\begin{rmk} \label{crucial} By construction of $\eta_X$ in Lemma \ref{fund}, $\omega_X^o$ is then the $\ell$-adic 
realization $r_{\ell, X}(\Sing (X))$ of the dg-category $\Sing(X)$, considered as 
a dg-category over $X$, and this is true without any regularity hypothesis on $X$. 
\end{rmk}

\begin{rmk} \label{fund-der} If $p: \mathcal{X} \to S$ is a proper lci map from a \emph{derived} scheme $\mathcal{X}$, we can still define a 2-periodic $\ell$-adic fundamental class map $\eta_{\mathcal{X}}$, 
as in Lemma \ref{fund}. This can be done by observing that the pushforward on G-theories along the inclusion of the truncation $\mathrm{t}_0 \mathcal{X} \to \mathcal{X}$ is an equivalence, and that $p$ being lci we have a natural inclusion $\Perf (\mathcal{X}) \to \Coh (\mathcal{X})$.
We further observe that in this case, while the $\s$-category $\D_c(\mathcal{X},\Ql(\beta))$ only depends
on the reduced subscheme $(\mathrm{t}_0\mathcal{X})_{red}$, and the same is true for the objects $\mathbb{Q}_{\ell, \mathcal{X}}(\beta)$ and
$\omega_{\mathcal{X}}(\beta)$, in contrast, the morphism $\eta_{\mathcal{X}}$ \emph{does} depend on 
the derived structure on $\mathcal{X}$, and thus it is not a purely topological invariant. 
\end{rmk}

Let us come back to $X$ a regular scheme, proper over $S$.
We have a canonical equivalence
in $\D_c(X_s,\Ql(\beta))$
\begin{equation}\label{final}
 \nu_X^\mathrm{I}(\beta)[1] \simeq \omega_{X_s}^o.
\end{equation}
This is a reformulation of the equivalence (\ref{prelim}), in view of Lemma \ref{fund}. 

\begin{rmk}\label{lcicase} When $X$ is not regular anymore, but still proper and lci\footnote{Note that a morphism of finite type between regular schemes is lci, since we can check that its relative cotangent complex has perfect amplitude in $[-1,0]$.} over $S$, 
there is nonetheless a natural morphism
$\nu_X^\mathrm{I}(\beta)[1] \longrightarrow \omega_{X_s}^o$, constructed as follows. 
Consider again the triangle (\ref{2})
$$\xymatrix{\nu_X^\mathrm{I} \ar[r] & \Ql \ar[r] & i_X^!(\Ql)[2].}$$
On $X$, we do have the
2-periodic $\ell$-adic fundamental class $\eta_X : \Ql(\beta) \longrightarrow \omega_X(\beta)$, and by taking its
$!$-pullback by $i_X^!$, we get a morphism 
$i_X^!(\Ql)(\beta) \longrightarrow \omega_{X_s}(\beta)$. This produces a 
sequence of morphisms
$$\xymatrix{\nu_X^\mathrm{I}(\beta) \ar[r] & \Ql(\beta) \ar[r] & i_X^!(\Ql)[2](\beta) = i_X^!(\Ql)
(\beta) \ar[r] & \omega_{X_s}(\beta).}$$
The resulting composite morphism $\Ql(\beta) \rightarrow \omega_{X_s}(\beta)$ is the 2-periodic $\ell$-adic  
fundamental class of $X_s$. Moreover, by construction, the composition
$\xymatrix{\nu_X^\mathrm{I}(\beta) \ar[r] & \Ql(\beta) \ar[r] &  \omega_{X_s}(\beta)}$
is canonically the zero map, and this induces the natural morphism
$$\alpha_X : \nu_X^{\mathrm{I}}(\beta)[1] \longrightarrow \omega_{X_s}^o$$ we were looking for.
Summing up, the morphism $\alpha_X$ always exists for any proper, lci scheme $X$ over $S$, and is an 
equivalence whenever $X$ is regular.
\end{rmk}

\subsection{A K\"unneth theorem for invariant vanishing cycles}

In this section, we consider two separated, finite type $S$-schemes $X$ and $Y$, 
such 
that $X_K$ and $Y_K$ are smooth over $K$, and both $X$ and $Y$ are regular and connected.
For simplicity\footnote{In the non-flat case, the fiber product of $X$ and $Y$ over $S$, to be considered below, should be replaced by the derived fiber product.} we also assume that $X$ and $Y$ are flat over $S$.
We set $Z:=X\times_S Y$, and consider this as a scheme over $S$. We have 
$Z_s \simeq X_s \times_s Y_s$, and the $\s$-category 
$D_c^\mathrm{I}(Z_s,\Ql)$ comes equipped with pull-back functors
$$p^* : D_c^\mathrm{I}(X_s,\Ql) \longrightarrow D_c^\mathrm{I}(Z_s,\Ql) \longleftarrow D_c^\mathrm{I}(Y_s,\Ql) :q^*.$$
By taking their tensor product, we get an external product functor
$$\boxtimes:=p^*(-)\otimes q^*(-) : D_c^\mathrm{I}(X_s,\Ql) \times D_c(Y_s,\Ql) \longrightarrow
D_c^\mathrm{I}(Z_s,\Ql).$$
For two objects $E \in  D_c^\mathrm{I}(X_s,\Ql)$ and $F \in  D_c^\mathrm{I}(Y_s,\Ql)$, we can define the K\"unneth
morphism in $D_c(Z_s,\Ql)$
$$\mathsf{k} : (E\boxtimes F)^\mathrm{I}[-1] \longrightarrow E^\mathrm{I} \boxtimes F^\mathrm{I}$$
as follows. Since Grothendieck duality on $Z_s$ is compatible with
external products, in order to define $\mathsf{k}$ it is enough to define its dual
$$\mathbb{D}(E^\mathrm{I}) \boxtimes \mathbb{D}(F^\mathrm{I}) \longrightarrow 
\mathbb{D}((E\boxtimes F)^\mathrm{I})[1].$$
By Lemma \ref{I-duality}, the datum of such a morphism is equivalent to that of a morphism
$$\mathbb{D}_I(E)^\mathrm{\mathrm{I}} \boxtimes \mathbb{D}_I(F)^\mathrm{I} \longrightarrow 
\mathbb{D}_I(E\boxtimes F)^\mathrm{I}.$$
We now define $\mathsf{k}$ as the map induced by the composite 

$$\xymatrix{\mathbb{D}_I(E)^\mathrm{I} \boxtimes \mathbb{D}_I(F)^\mathrm{I} \ar[r]^-{\mu_{(-)^\mathrm{I}}} & (\mathbb{D}_I(E) \boxtimes \mathbb{D}_I(F))^\mathrm{I}
\ar[r]^-{(\mu_{D_I})^\mathrm{I}} & \mathbb{D}_I(E\boxtimes F)^\mathrm{I}}$$
where $\mu_{(-)^\mathrm{I}}$ is the lax monoidal structure on $(-)^\mathrm{I}$, and $\mu_{D_I}$ the one on $D_I$\footnote{Note that $\mu_{D_I}$ is in fact an equivalence.}.

\begin{df}\label{conv}
With the above notations, the \emph{$\mathrm{I}$-invariant convolution}
of the two objects $E \in  D_c^\mathrm{I}(X_s,\Ql)$ and $F \in  D_c^\mathrm{I}(Y_s,\Ql)$ is defined to be the cone of the K\"unneth morphism, and 
denoted by $(E\circledast F)^\mathrm{I}$. By definition it sits in a triangle
$$\xymatrix{
(E\boxtimes F)^\mathrm{I}[-1] \ar[r]^-{\mathsf{k}} & E^\mathrm{I} \boxtimes F^\mathrm{I} \ar[r] & (E\circledast F)^\mathrm{I}.}$$
\end{df}

The main result of this section is the following proposition, relating 
the $\mathrm{I}$-invariant convolution of vanishing cycles on $X$ and $Y$ to the dualizing complex 
of $Z$. It can be also considered as a computation of the $\ell$-adic realization of 
the $\s$-category $\Sing(Z)$ of singularities of $Z$. \\
Note that, as $X$ and $Y$ are generically smooth over $S$, so is $Z$, and thus the 
2-periodic $\ell$-adic fundamental class map $\eta_Z : \Ql(\beta) \longrightarrow \omega_Z(\beta)$ of Lemma \ref{fund} is an equivalence 
over the generic fiber. Therefore, $\omega_Z^o$ is supported on $Z_s$, so that it can (and will)  be considered canonically as an object
in $\D_c(Z_s,\Ql)$. 

\begin{thm}\label{kunn}
With the above notations and assumptions, there is a canonical equivalence 
$$\omega_Z^o \simeq (\nu_X \circledast \nu_Y)^\mathrm{I}(\beta)$$ in  
$D_c(Z_s,\Ql(\beta))$.
\end{thm}

\noindent \textbf{Proof.} The proof of this theorem will combine various exact triangles together with an application of Gabber's K\"unneth formula for nearby cycles.

To start with, the 
vanishing cycles $\nu_Z$  of $Z$ sits in an exact triangle in $\D_c^\mathrm{I}(Z_s,\Ql)$
$$\xymatrix{ \Ql \ar[r] & \psi_Z \ar[r] & \nu_Z},$$ where
$\psi_Z=\psi_Z(\Ql)$ is the complex of nearby cycles of $Z$ over $S$. According to \cite[Lemma 5.1.1]{bb} or \cite[4.7]{illpadique}), we have
a natural equivalence in $\D_c^\mathrm{I}(Z_s,\Ql)$, induced by external product
$$\psi_Z \simeq \psi_X \boxtimes \psi_Y.$$
The object $\nu_Z$ then becomes the cone of the tensor product of the two morphisms in 
$\D_c^\mathrm{I}(Z_s,\Ql)$
$$\Ql \longrightarrow p^*(\psi_X) \qquad \Ql \longrightarrow q^*(\psi_Y)$$
where $p$ and $q$ are the two projections from $Z$ down to $X$ and $Y$, respectively.
Now, cones of tensor products are computed via the following well known lemma.

\begin{lem}\label{conetensor}
Let $\mathcal{C}$ be a stable symmetric monoidal $\s$-category, and 
$$u : x \rightarrow y \qquad v : x' \rightarrow y'$$
two morphisms. Let $C(u)$ be the cone of $u$, $C(v)$ be the cone of $v$, and $C(u\otimes v)$ the cone
of the tensor product $u\otimes v : x \otimes x' \rightarrow y \otimes y'$. 
Then, 
there exists a natural exact triangle
$$\xymatrix{C(u)\otimes x' \bigoplus  x\otimes C(v) \ar[r]&  C(u\otimes v) \ar[r] & 
C(u)\otimes C(v).}$$
\end{lem}

\noindent \emph{Proof of lemma.} Factor $u\otimes v$ as $\xymatrix{x\otimes x' \ar[r]^-{u\otimes id} & y\otimes x' \ar[r]^-{id \otimes v} & y\otimes y'}$, and apply the octahedral axiom to the triangles $$\xymatrix{x \otimes x' \ar[r]^-{u \otimes id} & y\otimes x' \ar[r]^-{f} & C(u)\otimes x' }$$
$$\xymatrix{y \otimes x' \ar[r]^-{id \otimes v} & y\otimes y' \ar[r] & y\otimes C(v) \ar[r]^-{d'}_-{[1]} & y\otimes x'[1] }, $$ to get a triangle 
$$\xymatrix{C(u) \otimes x' \ar[r] & C(u \otimes v) \ar[r] & y\otimes C(v) \ar[r]_-{[1]}^-{\theta} & C(u)\otimes x'[1]}$$ together with the compatibility
$\theta = f[1]\circ d'$. Now observe that $\theta \circ (u \otimes \mathrm{id}_{C(v)})=0$, and apply the octahedral axiom to the triangles
 $$\xymatrix{x \otimes C(v) \ar[r]^-{u \otimes id} & y\otimes C(v) \ar[r] & C(u)\otimes C(v) },$$
  $$\xymatrix{y \otimes C(v) \ar[r]^-{\theta} & C(u)\otimes x'[1] \ar[r]^-{f} & C(u\otimes v)[1]}$$ to conclude.
\hfill $\diamondsuit$ \\

Lemma \ref{conetensor} implies the existence of a natural exact triangle
in $\D_c^\mathrm{I}(Z_s,\Ql)$
$$\xymatrix{
\nu_X \boxplus \nu_Y \ar[r] & \nu_Z \ar[r] & \nu_X \boxtimes \nu_Y,}$$
which, by taking $\mathrm{I}$-invariants, yields an exact triangle in $\D_c(Z_s,\Ql)$
$$\mathbf{(T1)} \qquad \xymatrix{
\nu_X^\mathrm{I} \boxplus \nu_Y^\mathrm{I} \ar[r] & \nu_Z^\mathrm{I} \ar[r] & (\nu_X \boxtimes \nu_Y)^\mathrm{I}.}$$

\begin{lem}\label{dualizfiber} Let $k$ be an algebraically closed field, $s:= \mathrm{Spec} \, k$, $ p_X: X \to s$, $p_Y: Y \to s $ be proper morphisms of schemes, and $p_1: Z:= X \times_s Y \to X$, $p_2: Z:= X \times_s Y \to Y$ the natural projections.  If $\omega_Z$, $\omega_X$, $\omega_Y$ denote the $\Ql$-adic dualizing complexes of $Z$, $X$, and $Y$, respectively, there is a canonical equivalence  $$a: p_1^*\omega_X \otimes p_2^*\omega_Y \longrightarrow \omega_Z. $$
\end{lem}

\noindent \emph{Proof of lemma.} We first exhibit the map $a$. We denote simply by $\underline{\mathrm{Hom}}_T(-, -)$ the derived internal hom in $\D_{c}(T,\Ql)$ (so that $D:= \underline{\mathrm{Hom}}_Z(-, \omega_Z)$ is the $\Ql$-adic duality on $Z$). By adjunction, giving a map $a$ is the same thing as giving a map $ \omega_X \to (p_1)_*\underline{\mathrm{Hom}}_Z(p_2^* \omega_Y, p_2^! \omega_Y )$. Since $p_X$ is proper, by \cite[Exp XVIII, 3.1.12]{sga4III}, we have a canonical equivalence $$(p_1)_*\underline{\mathrm{Hom}}_Z(p_2^* \omega_Y, p_2^! \omega_Y ) \longrightarrow p_X^! (p_Y)_* \underline{\mathrm{Hom}}_Y(\omega_Y, \omega_Y).$$ Therefore, we are left to define a map $$\omega_X \simeq p_X^! \Ql \to p_X^! (p_Y)_* \underline{\mathrm{Hom}}_Y(\omega_Y, \omega_Y),$$ and we take $p_X^!(\alpha)$ for this map, where $\alpha$ is the adjoint to the canonical map $\Ql \simeq p_Y^*\Ql \to \underline{\mathrm{Hom}}_Y(\omega_Y, \omega_Y)$.  One can then prove that $a$ is an equivalence, by checking it stalkwise. \hfill $\Box$ \\

By Lemma \ref{dualizfiber}, the dualizing
complex $\omega_{Z_s}= (Z_s\to s)^! \Ql$ of $Z_s \simeq X_s \times_s Y_s$ is canonically equivalent to 
$\omega_{X_s} \boxtimes \omega_{Y_s}$, and, through this equivalence, the virtual fundamental 
class of $Z_s$ 
$$\eta_{Z_s} : \Ql(\beta) \longrightarrow \omega_{Z_s}(\beta) \simeq 
\omega_{X_s} \boxtimes \omega_{Y_s} (\beta)$$
is simply given by the external tensor product of the virtual fundamental classes of 
$X_s$ and $Y_s$. By Lemma \ref{conetensor} we thus get 
a second exact triangle in $\D_c(Z_s,\Ql(\beta))$
$$\mathbf{(T2)} \qquad \xymatrix{
\omega_{X_s}^o \boxplus \omega_{Y_s}^o \ar[r] & \omega_{Z_s}^o \ar[r] & 
\omega_{X_s}^o \boxtimes \omega_{Y_s}^o.}$$
There is a morphism from the triangle (T1) to the triangle (T2) which is defined 
using the natural morphism
$$\alpha_Z : \nu_Z^\mathrm{I}(\beta)[1] \longrightarrow \omega_{Z_s}^o$$
introduced in Remark \ref{lcicase}.  In fact, $Z$ is proper and lci over $S$ (since $X/S$ is flat and lci\footnote{Again, note that a morphism of finite type between regular schemes is lci, since we can check that its relative cotangent complex has perfect amplitude in $[-1,0]$.}, and being lci is stable under flat base change and composition), and 
the map $\alpha_Z$ is defined for any proper, lci scheme $Z$ over $S$,  
being an equivalence when $Z$ is regular with smooth
generic fiber. Using the compatible maps
$\alpha_X$, $\alpha_Y$ and $\alpha_Z$ we get a commutative square
$$\xymatrix{
\nu_X^\mathrm{I}(\beta)[1] \boxplus \nu_Y^\mathrm{I}(\beta)[1] \ar[d]_-{\alpha_Z \oplus \alpha_Y}\ar[r] & 
\nu_Z^\mathrm{I}(\beta)[1] 
\ar[d]^-{\alpha_Z} \\
\omega_{X_s}^o \boxplus \omega_{Y_s}^o \ar[r] & \omega_{Z_s}^o.}$$
This produces a morphism from triangle $(T1)$ (tensored by $\Ql(\beta)[1]$) to 
triangle $(T2)$ 
\begin{equation}\label{mapp}\xymatrix{
\nu_X^\mathrm{I}(\beta)[1] \boxplus \nu_Y^\mathrm{I}(\beta)[1] \ar[r] \ar[d]  & 
\nu_Z^\mathrm{I}(\beta)[1] \ar[r]  \ar[d]^-{\alpha_Z} & (\nu_X \boxtimes \nu_Y)^\mathrm{I}(\beta)[1]
\ar[d] \\
\omega_{X_s}^o \boxplus \omega_{Y_s}^o \ar[r] & \omega_{Z_s}^o \ar[r] & 
\omega_{X_s}^o \boxtimes \omega_{Y_s}^o.}
\end{equation}
Since $X$ and $Y$ are regular
with smooth generic fibers, the maps $\alpha_X$ and $\alpha_Y$ are equivalences, therefore the leftmost vertical morphism is also an equivalence.
Thus the right hand square is a cartesian square. 

Now, the rightmost vertical morphism can be written, again using the equivalences $\alpha_X$ and $\alpha_Y$, as
$$(\nu_X \boxtimes \nu_Y)^\mathrm{I}(\beta)[1] \longrightarrow (\nu_X^\mathrm{I}(\beta)[1]) 
\boxtimes_{\Ql(\beta)} (\nu_Y^\mathrm{I}(\beta)[1]) \simeq (\nu_X^\mathrm{I} \boxtimes \nu_Y^\mathrm{I})[2](\beta)$$
This morphism is the K\"unneth map $\mathsf{k}$ of Definition \ref{conv} tensored by $\Ql[2](\beta)\simeq \Ql(\beta)$, 
and thus its cone is
$(\nu_X \circledast \nu_Y)^\mathrm{I}(\beta)$. In order to finish the proof of the proposition it then  
remains to show that the cone of the middle vertical morphism in $(\ref{mapp})$
$$\alpha_Z: \nu_Z^\mathrm{I}(\beta)[1] \longrightarrow \omega_{Z_s}^o$$
can be canonically identified with $\omega_Z^o$.

For this, we remind the exact triangle $(\ref{2})$ tensored by $\Ql (\beta)$
$$\xymatrix{
\Ql(\beta) \ar[r] & i_Z^!(\Ql(\beta)) \ar[r] & \nu_Z^\mathrm{I}(\beta)[1].}$$
Using the fundamental class $\eta_Z : \Ql(\beta) \rightarrow \omega_Z(\beta)$, we get 
a morphism of triangles
$$\xymatrix{
 \Ql(\beta) \ar[r] \ar[d]_-{id} & i_Z^!(\Ql(\beta)) \ar[r] \ar[d]^-{i_Z^!(\eta_Z)} & 
 \nu_Z^\mathrm{I}(\beta)[1] \ar[d] \\
 \Ql(\beta) \ar[r] & i_Z^!(\omega_Z(\beta)) \ar[r] & \omega_{Z_s}^o.}$$
The right hand square is thus cartesian, so that the cone of the vertical 
morphism on the right is canonically identified with the cone of the vertical 
morphism in the middle. By definition, this cone is
$i_Z^!(\omega_Z^o)$. Since $Z_{K}$ is smooth over $K$, the $\ell$-adic complex $\omega_Z^o$ is supported
on $Z_s$, and thus $i_Z^!(\omega_Z^o)$ is canonically equivalent to $\omega_Z^o$, and we conclude.
\hfill $\Box$ \\

\begin{cor}\label{ckunn}
We keep the same notations and assumptions as in Proposition \ref{kunn}, and we further assume 
one of the following conditions:

\begin{enumerate}

\item the $\mathrm{I}$-action on $\nu_X$ and on $\nu_Y$ is tame, or

\item the reduced scheme $(X_s)_{red}$ is smooth over $k$. 

\end{enumerate}

Then, there is a canonical equivalence
$$\omega_Z^o \simeq (\nu_X \boxtimes \nu_Y)^\mathrm{I}(\beta)$$
 in $\D_c(Z_s,\Ql(\beta))$.
\end{cor}

\noindent\textbf{Proof.} It is enough to prove that under any one of the two assumptions, 
 $(\nu_X \circledast \nu_Y)^\mathrm{I}$ is canonically equivalent to 
$(\nu_X \boxtimes \nu_Y)^\mathrm{I}$. If the scheme $(X_s)_{red}$ is smooth over $k$, then 
we have $\nu_X^\mathrm{I}(\beta)=0$. Indeed, triangle (\ref{2}) can be then re-written 
$$\xymatrix{
\nu_X^\mathrm{I} \ar[r] & \Ql \ar[r] & \Ql[2n+2]}$$ where $n$ is the dimension of $X_s$.
By tensoring by $\Ql (\beta)$, we get a triangle $$\xymatrix{
\nu_X^\mathrm{I} (\beta) \ar[r] & \Ql(\beta) \ar[r]^-{b} & \Ql[2n+2](\beta) \simeq \Ql (\beta)}$$
where $b$ is an equivalence. Therefore, $\nu_X^\mathrm{I} (\beta) =0$, and, 
by definition of $\mathrm{I}$-invariant convolution, this implies that 
$(\nu_X \circledast \nu_Y)^\mathrm{I} \simeq (\nu_X \boxtimes \nu_Y)^\mathrm{I}$.

Assume now that the action of $\mathrm{I}$ on $\nu_X$ and $\nu_Y$ is tame. This means that 
the action of $\mathrm{I}$ factors through the natural quotient $\mathrm{I} \longrightarrow \mathrm{I}_t$, where 
$\mathrm{I}_t$ is the tame inertia group, which is canonically isomorphic
to $\hat{\mathbb{Z}}'$, the prime-to-$p$ part
of the profinite completion of $\mathbb{Z}$. As we have chosen 
a topological generator $T$ of $\mathrm{I}_t$ (see Section \ref{trivtate}), 
the actions of $\mathrm{I}$ are then 
completely caracterized by the automorphisms $T$ on $\nu_X$ and $\nu_Y$. 
Moreover, $\nu_X^\mathrm{I}$ is then naturally equivalent to the homotopy fiber of
$(1-T) : \nu_X \rightarrow \nu_X$, and similarly for $\nu_Y^\mathrm{I}$. From this it is easy to see that 
the K\"unneth map
$$(\nu_X \boxtimes \nu_Y)^\mathrm{I}[-1] \longrightarrow \nu_X^\mathrm{I} \boxtimes \nu_Y^\mathrm{I}$$
fits in an exact triangle
$$\xymatrix{
(\nu_X \boxtimes \nu_Y)^\mathrm{I}[-1] \ar[r] & 
\nu_X^\mathrm{I} \boxtimes \nu_Y^\mathrm{I} \ar[r] & (\nu_X \boxtimes \nu_Y)^\mathrm{I}}$$
where the second morphism is induced by the lax monoidal structure on $(-)^\mathrm{I}$.
We conclude that there is a natural equivalence
$(\nu_X \circledast \nu_Y)^\mathrm{I} \simeq (\nu_X \boxtimes \nu_Y)^\mathrm{I}$.
\hfill $\Box$ \\

\section{K\"unneth formula for dg-categories of singularities}

\subsection{The monoidal dg-category $\B$ and its action}

We keep our standing assumptions: $A$ is an excellent strictly henselian dvr with perfect residue field $k$ and
fraction field $K$. We let $S= \mathrm{Spec}, A$ and $s=\mathrm{Spec}\, k$, as usual, and choose an uniformizer
$\pi$ of $A$.

We let $G:=s\times^h_S s$ (derived fiber product), considered as a derived scheme over $S$. The derived scheme $G$ 
has a canonical structure of groupoid in derived schemes acting on $s$. 
The composition in the groupoid $G$ induces a convolution monoidal structure on 
the dg-category of coherent complexes on $G$
$$\odot : \Coh(G) \otimes_A \Coh(G) \longrightarrow \Coh(G).$$
More explicitly, we have a map of derived schemes
$$\xymatrix{G\times_s G \ar[r]^-{q} & G,}$$
defined as the projection on the first and third components
$s\times_S s \times_S s \rightarrow s\times_S s$. We then  
define $\odot$ by the formula
$$E \odot F := q_*(E\boxtimes_s F)$$
for two coherent complexes $E$ and $F$ on $G$. More generally, 
if $X \longrightarrow S$ is any scheme, with special fiber $X_s$ (possibly a derived scheme, by taking the derived fiber at $s$), 
the groupoid $G$ acts naturally on $X_s$ via the natural projection
$$q_X : G\times_s X_s \simeq (s\times_S s)\times_s (s\times_S X) \simeq s \times_S X_s \longrightarrow
X_s.$$
This defines an external action
$$\odot : \Coh(G) \otimes_A \Coh(X_s) \longrightarrow \Coh(X_s)$$
by 
$E\odot M := (q_X)_*(E\boxtimes_s M)$.

The homotopy coherences issues for the above $\odot$-structures can be handled using the fact that 
the construction $Y \mapsto \Coh(Y)$ is in fact a symmetric lax monoidal $\s$-functor
from a certain $\s$-category of correspondences between derived schemes
to the $\s$-category of dg-categories over $A$. As a result, $\Coh(G)$ is endowed with a natural structure of a monoid in 
the symmetric monoidal $\s$-category 
$\dgcat_A$, and that, for any $X/S$, $\Coh(X_s)$ is naturally a module over 
$\Coh(G)$ in $\dgcat_A$. However, for our purposes it will be easier and more efficient
to provide explicit models for both $\Coh(G)$ and its action on 
$\Coh(X_s)$. This will be done locally in the Zariski topology 
in a similar spirit to \cite[Section 2]{brtv}; the global construction will then be obtained by 
a rather straightforward gluing procedure. \\

\subsubsection{The monoidal dg-categories $\B^+$ and $\B$} Let 
$K_A$ be the Koszul commutative $A$-dg-algebra of $A$ with respect to $\pi$
$$K_A : \xymatrix{A \ar[r]^-{\pi} & A}$$
sitting in degrees $[-1,0]$.
The canonical generator of $K_A$ in degree $-1$ will be denoted by $h$. 
In the same way, we define the commutative $A$-dg-algebra
$$K^2_A:=K_A \otimes_A K_A$$
which is the Koszul dg-algebra of $A$ with respect to the sequence $(\pi,\pi)$. 
As a commutative graded $A$-algebra, $K^2_A$ is 
$Sym_A(A^2[1])$, and it is endowed 
with the unique multiplicative differential sending the two generators $h$ and $h'$ in degree $-1$ 
to $\pi$ (and $hh'$ to $\pi \cdot h' - \pi \cdot h$).

Moreover, $K^2_A$ has a canonical structure of \emph{Hopf algebroid} over $K_A$, in which the source and target map are the two natural inclusions
$K_A \longrightarrow K^2_A$, whereas the unit is given by the multiplication $K^2_A \rightarrow K_A$. 
The composition (or coproduct) in this Hopf algebroid structure is given by 
$$\Delta:= \mathrm{id}\otimes 1 \otimes \mathrm{id} : K_A \otimes_A K_A=K^2_A \longrightarrow K^2_A \otimes_{K_A}K^2_A = K_A \otimes_A K_A \otimes_A K_A .$$
Finally, the antipode is the automorphism 
of $K^2_A$ exchanging the two factors $K_A$. This structure of Hopf algebroid
endows the dg-category $\mathsf{Mod}(K^2_A)$, of $K^2(A)$-dg-modules, 
with a unital and associative monoidal structure $\odot$. It is explicitely given 
for two object $E$ and $F$, by the formula
$$E\odot F := E\otimes_{K_A}F$$
where the $K_A \otimes_A K_A \otimes_A K_A$-module on the rhs is considered
as a $K^2_A$ module via the map $\Delta$. The unit of this monoidal structure 
is the object $K_A$, viewed as a $K^2_A$-module by the multiplication $K^2_A \rightarrow K_A$. 
It is not hard to see that $\odot$ preserves 
cofibrant $K^2_A$-dg-modules; more generally it makes $\mathsf{Mod}(K^2_A)$ into
a monoidal model category in the sense of \cite[Ch. 4]{ho}. Note however that 
the unit $K_A$ is \emph{not} cofibrant in this model structure.

\begin{df}
The \emph{monoidal dg-category $\B_{str}^+$} is defined to be
$\Mod^c(K^2_A)$, the dg-category of all
cofibrant dg-modules over $K^2_A$ which are perfect over $A$,
together with the unit object $K_A$. It is endowed
with the monoidal structure $\odot$ described above. 
\end{df}

By Appendix A, the localization of the dg-category $\B_{str}^+$ along all
quasi-isomorphisms, defines a monoidal dg-category.

\begin{df}
The \emph{monoidal dg-category $\B^+$} is defined to be the localization 
$$W_{\textrm{eq}}^{-1}(\B^+_{str}),$$
where $W_{\textrm{eq}}$ is the set of quasi-isomorphisms. It is naturally a unital and associative
monoid in the symmetric monoidal $\s$-category $\dgcat_A$.
\end{df}

\begin{rmk}\label{thisisamodel}Note that $\B^+$ defined above is a model for $\Coh(G)$, for our derived groupoid
$G=s\times_S s$ above. Indeed, the commutative dga $K^2_A$ is 
quasi-isomorphic to the normalization of the simplicial  algebra $k\otimes^{\mathbb{L}}_A k$. 
Now, since $G \to S$ is a closed immersion,
 a quasi-coherent complex $E$ on $G$ is coherent iff its direct image on $S$ is coherent, hence perfect, $S$ being regular.
In particular, we have that $\Coh(G)$ is equivalent to the dg-category of
all cofibrant $K^2_A$-dg-modules which are perfect over $A$ 
The latter dg-category is also naturally equivalent to the
localization of $\B^+_{str}$ along quasi-isomorphisms\footnote{We leave to the
reader to check that adding the unit objet $K_A$ does not change the
localization.} 
\end{rmk}

We now introduce the weak monoidal dg-category $\B$, defined as a further
localization of $\B^+$. This will be our main ``base monoid'' for the module dg-categories we will be interested in.

\begin{df}
The \emph{weak monoidal dg-category $\B$} is defined to be the localization
$$\B := L_{W}(\B^+),$$
where $W$ is the set of morphisms in $\B^+$ 
whose cone is perfect as a $K^2_A$-dg-module.
\end{df}

As $\B^+$ it itself defined as a localization of $\B^+_{str}$,
if $W_{\textrm{pe}}$ denotes the morphisms in  $\B^+_{str}$  
whose cones are perfect over $K^2_A$, we have $W_{\textrm{eq}} \subseteq W_{\textrm{pe}}$, 
then $\B$ can also be realized as localization of $\B^+_{str}$ directly, as 
$$\B \simeq L_{W_{\textrm{pe}}}\B^+_{str}.$$
Since the monoidal structure $\odot$ is compatible 
 with $W_{\textrm{pe}}$, this presentation implies that $\B$ comes equipped with a natural 
structure of an associative and unital
monoid in $\dgcat_A$. Note, moreover, that $\B$ comes equipped with a natural morphism of monoids
$$\B^+ \longrightarrow \B$$ given by the localization map. \\

\subsubsection{The local actions}\label{localact}

Let now $X=\mathrm{Spec}\, R$ be a regular scheme flat over $A$. As done above for the monoid structure on $\B^+$, we will define a \emph{strict
model} for $\Coh(X_s)$, together with a \emph{strict model} for the $\Coh(G)$-action on $\Coh(X_s)$. In order to do this, 
let $K_R$ be the Koszul dg-algebra of $R$ with respect to $\pi$, which 
comes equipped with a natural map $K_A \longrightarrow K_R$ of cdga's over $A$. We consider
$\Mod^c(K_R)$, the dg-category of all cofibrant $K_R$-dg-modules which are
perfect as $R$-modules (note that $R$ is regular, and see Remark \ref{thisisamodel}). 
The same argument as in Remark \ref{thisisamodel} then shows that this dg-category is naturally equivalent to $\Coh(X_s)$.
Moreover, $\Mod^c(K_R)$ has a structure of 
a $\B^{+}_{str}$-module dg-category defined as follows. For $E \in \B_{str}^+$, and  
$M \in \Mod^c(K_R)$, we can define
$$E\odot M := E\otimes_{K_A}M,$$
where, in the rhs, we have used the ``right'' $K_A$-dg-module structure on $E$, i.e. the one induced by the composition $$\xymatrix{K_A \ar[r]^-{\sim} & A\otimes_A K_A \ar[r]^-{\mathrm{id}\otimes u} & K_A\otimes_A K_A },$$  $u: A \to K_A$ being the canonical map. As $E$ is either the unit or
it is cofibrant over $K_A^2$ (and thus cofibrant over $K_A$), $E\otimes_{K_A}M$ is again
a cofibrant $K_R$-module, and again perfect over $R$, i.e. $E\odot M \in \Mod^c(K_R)$. By localization along quasi-isomorphisms, 
(see Proposition \ref{detailslocalactions} for details) 
we obtain that $\Coh(X_s)$
carries a natural $\B^+$-module structure as 
an object in the symmetric monoidal $\s$-category $\dgcat_A$. \\

We now apply a similar argument in order to define a \emph{$\B$-action on $\Sing(X_s)$}.
Let again $X=Spec\, R$ be a regular scheme over $A$, and consider 
$\Mod^c(K_R)$ as a $\B^+_{str}$-module dg-category as above. Let $W_{R, \textrm{pe}}$ be the set 
of morphisms in $\Mod^c(K_R)$ whose cones are perfect dg-modules over $K_R$. By localization we then get 
a $\B$-module structure on $L_{W_{R, \textrm{pe}}}\Mod^c(K_R)$. Note that the localization $L_{W_{R, \textrm{pe}}}\Mod^c(K_R)$
 is a model for the dg-category $\Sing(X_s)$, which therefore comes equipped with 
the structure of a $\B$-module in $\dgcat_A$. \\

We gather the details of above constructions in the following

\begin{prop}\label{detailslocalactions} Let $X= \mathrm{Spec}\, R$ be a regular scheme, flat over $S= \mathrm{Spec} \, A$, and $X_s$ its special fiber.
Then there is a canonical $\B^+$-module structure (resp., $\B$-module structure) on $\Coh (X_s)$ (resp., on $\Sing (X_s)$), inside $\dgcat_A$.
\end{prop}

\noindent \textbf{Proof.} This is an easy application of the localization results presented in Appendix A.\\ We first treat the case of $\B^+$ and $T:= \Coh (X_s)$. If $T^{\textrm{str}}:= \Mod^c(X_r)$ and $W_{T, \textrm{eq}}$ denotes the quasi-isomorphisms in $T^{\textrm{str}}$, we have $W_{T, \textrm{eq}}^{-1}T^{\mathrm{str}} \simeq T$ in $\dgcat_A$. Analogously, $W_{\mathrm{eq}}^{-1} \B_{\mathrm{str}}^+ \simeq \B^+$ in $\dgcat_A$. In order to apply the localization result of Appendix A, we need to prove that the tensor product $\odot: \B^+_{\str} \otimes_A T^{\str}\to T^{\str}$ (defined in \ref{localact}) sends $W_{\textrm{eq}}\otimes id \cup id \otimes W_{T, \textrm{eq}}$ to $W_{T, \textrm{eq}}$. If $L, L' \in \B^+_{\str}$, $E, E' \in T^{\str}$, and $w': L\to L'$, $w: E \to E'$ are quasi-isomorphisms, then  $w'\odot id_E$ is again a quasi-isomorphism (because $L$, and $L'$ are cofibrant over $K_A^2$, hence over $K_A$, and thus $w'$ is in fact a homotopy equivalence), and the same is true for $id_L \odot w$ (since $L$ is cofibrant over $K_A$). Therefore, $\odot$ does send $W_{\textrm{eq}}\otimes id \cup id \otimes W_{T, \textrm{eq}}$ to $W_{T, \textrm{eq}}$, and there is an induced canonical map $(W_{\textrm{eq}}\otimes id \cup id \otimes W_{T, \textrm{eq}})^{-1} \B^+_{\str} \otimes_A T^{\str} \to W_{T, \textrm{eq}}^{-1} T^{\str} \simeq \Coh (X_s)$. By composing this with the natural equivalence  $(W_{\textrm{eq}}\otimes id \cup id \otimes W_{T, \textrm{eq}})^{-1} \B^+_{\str} \otimes_A T^{\str} \to W_{\mathrm{eq}}^{-1} \B_{\mathrm{str}}^+ \otimes_A W_{T, \textrm{eq}}^{-1}T^{\mathrm{str}}$ (Appendix A), we finally get our $\B^+$-module structure on $\Coh (X_s)$ inside $\dgcat_A$. \\We now treat the case of $\B$ and $T= \Sing (X_s)$. Here, we consider the pairs $(T^{\str}= \Mod^c(X_r), W_{T})$ where $W_T$  are the maps in $T^{\str}$ whose cones are perfect over $K_R$, and $(\B^*_{\str}, W)$, where $W$ are the maps in $\B^+_{\str}$ whose cones are perfect over $K^2_A$. We have $W^{-1}\B^*_{\str} \simeq \B$, and $W_T^{-1} T^{\str} \simeq \Sing (X_s)$, and we need to prove that both $W\odot id$ and $id \odot W_T$ are contained in $W_T$. Let $u: L\to L' \in W$ and $v: E \to E' \in W_T$, and $C(-)$ denote the cone construction. We have $C(id_L \odot v) \simeq L\otimes_{K_A} C(v)$ and $C(u \odot id_E) \simeq C(u) \otimes_{K_A} E$. By hypothesis, $L$ is perfect over $A$ hence over $K_A$ (since $K_A$ is perfect over $A$), and $C(v)$ is perfect over $K_A$, since $X \to S$ is lci (as a map of finite type between regular schemes), and thus $K_A \to K_R$ is derived lci (recall that $X/S$ is flat so that $X_s$ is also the derived fiber), and pushforward along a lci map preserves perfect complexes. So,  $C(id_L \odot v) \in W_T$. On the other hand, the ``right-hand'' map $K_A \to K^2_A $ (with respect to which $L$ and $L'$ are viewed as $K_A$-dg-modules in the definition of $\odot$) is derived lci, hence $C(u)$ is perfect over $K_A$, being perfect over $K^2_A$ by hypothesis. Moreover,  since $s\to S$ is a closed immersion, $E$ is perfect over $K_A$ iff $(X\to S)_* E$ is perfect (= coherent, $S$ being regular) over $S$; but $(X_s \to X)_* E$ is perfect by hypothesis, and pushforward along $X\to S$ preserves perfect complexes, since $X/S$ is lci. Therefore  $C(u \odot id_E) \in W_T$, and we deduce that $\odot$ does send $W\otimes id \cup id \otimes W_T$ to $W_T$. This gives us an induced canonical map $(W\otimes id \cup id \otimes W_T)^{-1} \B^+_{\str} \otimes_A T^{\str} \to W_T^{-1} T^{\str} \simeq \Sing (X_s)$. By composing this with the natural equivalence  $(W\otimes id \cup id \otimes W_T)^{-1} \B^+_{\str} \otimes_A T^{\str} \to W^{-1} \B_{\mathrm{str}}^+ \otimes_A W_T^{-1}T^{\mathrm{str}}$ (Appendix A), we finally get our $\B$-module structure on $\Sing (X_s)$ inside $\dgcat_A$. 

\hfill $\Box$ \\

\begin{lem}\label{lembrtv}
The natural morphism
$$ \Coh(X_s)^o \otimes_{\B^+}\B \longrightarrow \Sing(X_s) $$
is an equivalence of $\B$-modules.
\end{lem}

\noindent \textbf{Proof.} This is a reformulation of \cite[Proposition 2.43]{brtv}. \hfill $\Box$ \\

\subsubsection{The global actions} We now let $X$ be a regular scheme, not necessarily affine anymore. 
We have by Zariski descent
$$\Coh(X_s) \simeq \lim_{Spec\, R \,\subset \, X}\Mod^c(K_R),$$
where the limit is taken over all affine opens $Spec\, R$ of $X$.
The right hand side of the above equivalence is a limit of dg-categories
underlying $\B^+$-module structures (Proposition \ref{detailslocalactions}). As the forgetful functor
from $\B^+$-modules to dg-categories reflects limits, this endows
$\Coh(X_s)$ with a unique structure of $\B^+$-module. 

In the same way, we have Zariski descent for $\Sing(X_s)$ in the sense that 
$$\Sing(X_s) \simeq \lim_{Spec\, R \subset X}L_{W_{R, \textrm{pe}}}(\Mod^c(K_R)).$$
The right hand side of the above equivalence is a limit of dg-categories
underlying $\B$-module structures (Proposition \ref{detailslocalactions}). As the forgetful functor
from $\B$-modules to dg-categories reflects limits, this makes 
the dg-category $\Sing(X_s)$ into a $\B$-module in a natural way. \\

An important property of these $\B^+$ and $\B$-module structures is given in the following proposition.

\begin{prop}\label{cotens}
Let $X$ be a flat regular scheme over $S$.
\begin{enumerate}
\item 
The $\B^+$-module structure on $\Coh(X_s)$ is cotensored. 
\item The $\B$-module structure on $\Sing(X_s)$ is cotensored. 
\end{enumerate}
\end{prop}

\noindent \textbf{Proof.} This follows from Remark \ref{little} since  
the monoidal triangulated dg-categories $\B^+$ and $\B$ are both generated by their unit objects.
\hfill $\Box$ \\

\subsection{K\"unneth formula for dg-categories of singularities}

In the previous section we have seen that,  
for any regular scheme $X$ over $S$, the dg-category $\Sing(X_s)$ are equipped
with a natural $\B$-module structure. In this section we compute 
tensor products of dg-categories of singularities over $\B$.  \\

From a general point of view, let $T \in \dgcat_A$ be a $\B$-module, and assume
that $T$ is also \emph{co-tensored} over $\B$ (Def. \ref{cotens}). Then $T^{o}$ has a natural structure
of a $\B^{\otimes-op}$-module given by co-tensorisation. By Proposition \ref{cotens}, we may take $T=\Sing(X_s)$,
so that $T^o$ is a $\B^{\otimes-op}$-module. 
In particular, if we have another regular scheme $Y$, we are entitled to consider the tensor product
$$\Sing(X_s)^o \otimes_{\B}\Sing(Y_s),$$
which is a well defined object in $\dgcat_A$. We further assume, for simplicity, that $X$ and $Y$ are
also flat over $S$.
The main result of this section is the following proposition, which is a
categorical counterpart of our K\"unneth formula for vanishing cycles (Proposition \ref{kunn}).

\begin{thm}\label{dg-kunn}
Let $X$ and $Y$ be two regular schemes, flat over $S$, with smooth 
generic geometric fibers. 
There is a natural equivalence in $\dgcat_A$
$$\Sing(X_s)^o \otimes_{\B} \Sing(Y_s) \simeq \Sing(X \times_S Y).$$
\end{thm}
 
\noindent \textbf{Proof.} 
As $\B$ is a localization of $\B^+$, the 
natural $\s$-functor on enriched dg-categories (which is lax symmetric monoidal)
$$\dgcat_\B \longrightarrow \dgcat_{\B^+}$$
is fully faithful. Moreover, it commutes with tensor products in the following
sense: for two objects $T$ and $T'$ in $\dgcat_\B$, there is a natural 
equivalence
$T^o\otimes_\B T' \simeq T^o\otimes_{\B^+}T'.$
Therefore, in order to prove the theorem it enough to construct an equivalence
of dg-categories over $A$
$$\Sing(X_s)^o \otimes_{\B^+} \Sing(Y_s) \simeq \Sing(X \times_S Y).$$

Now, we claim that the result is local on 
$Z:=X\times_S Y$.  Indeed, we have two prestacks of dg-categories on 
the small Zariski site $Z_{zar}$: 
$$U\times_S V \mapsto \Sing(U_s)^o \otimes_{\B} \Sing(V_s) \qquad
U\times_S V \mapsto \Sing(U\times_S V)$$
for any two affine opens $U\subset X$ and $V\subset Y$. These two 
prestacks are stacks of dg-categories, and thus we have
equivalences in $\dgcat_A$
\begin{equation}\label{produ}\Sing(X_s)^o \otimes_{\B} \Sing(Y_s) \simeq \lim_{U\subset X, V \subset Y}
\Sing(U_s)^o \otimes_{\B} \Sing(V_s)\end{equation}
\begin{equation}\label{singzar}\Sing(X\times_S Y) \simeq \lim_{U\subset X, V \subset Y}
\Sing(U\times_S V).\end{equation}
The stack property (\ref{singzar}) is proved in \cite[2.3]{brtv} (where it is moreover 
shown that this is a stack for the h-topology). The stack property (\ref{produ}) is then a consequence of the same descent argument for dg-categories of singularities. Indeed, we have the following lemma. 

\begin{lem}\label{lem1}
Let $Z$ be an $S$-scheme, and $F$ be a stack
of  $\mathcal{O}_Z$-linear dg-categories. Assume that $F$ 
is a $\B^{\otimes-op}$-module stack\footnote{I.e., $F$ is a stack on $Z_{\textrm{Zar}}$ with values in the $\infty$-category of $\B^{\otimes-op}$-modules in $\dgcat_A$.}, and let $T_0$ be a $\B$-module dg-category.
Then, the prestack $F \otimes_{\B}T_0$ of dg-categories of $Z_{Zar}$, sending 
$W\subset Z$ to $F(W)\otimes_{\B}T_0$ is a stack.
\end{lem}

\noindent \textbf{Proof of Lemma \ref{lem1}.} This is an application of the main result of 
\cite{to2}. We denote, as usual, by 
$\widehat{T}:=\rch(T,\widehat{A})$ the (non-small) dg-category of all $T^o$-dg-modules. 
By the main result of \cite{to2}, the dg-category  
$$\lim_{W\subset Z}
\widehat{(F(W)\otimes_{\B}T_0)}$$
is compactly generated, and its dg-category of compact objects is 
equivalent in $\dgcat_A$ to $\lim_{W\subset Z} F(W)\otimes_{\B}T_0$. 
Moreover, we have
$$\widehat{(F(W)\otimes_{\B}T_0)} \simeq \widehat{F(W)} \widehat{\otimes}_{\widehat{\B}}
\widehat{T_0},$$
where $\widehat{\otimes}$ is the symmetric monoidal structure
on presentable dg-categories. As $\widehat{\otimes}$ is rigid
when restricted to compactly generated dg-categories, we have that 
$\widehat{\otimes}_{\widehat{\B}}$ distributes over limits on both factors, and we have
$$\widehat{(F(Z)\otimes_{\B}T_0)}
\simeq \lim_{W\subset Z} \widehat{F(W)} \widehat{\otimes}_{\widehat{\B}}\widehat{T_0}.$$
Passing to the sub-dg-categories of compact objects, we find that 
$$F(Z) \otimes_{\B}T_0 \simeq \lim_{W\subset Z}F(W) \otimes_{\B}T_0$$
which is the statement of the lemma. \hfill $\Diamond$ \\

Lemma \ref{lem1} immediately implies the stack property (\ref{produ}).\\  
We are thus reduced to the case where $X$ and $Y$
are both affine, and we have to produce an equivalence
$$\Sing(X_s)^o \otimes_{\B} \Sing(Y_s) \simeq \Sing(X \times_S Y)$$
that is compatible with Zariski localization on both $X$ and $Y$. \\

In this case, we start by the following.

\begin{lem}\label{lem2}
There is a natural equivalence of dg-categories over $A$
$$\Coh(X_s)^o \otimes_{\B^+} \Coh(Y_s) \simeq \Coh_{Z_s}(Z),$$
where $Z:=X\times_S Y$, and the right hand side is the dg-category of
coherent complexes on $Z$ with cohomology supported on the special fiber $Z_s$. 
This equivalence is furthermore functorial in $X$ and $Y$. 
\end{lem}

\noindent \textbf{Proof of lemma \ref{lem2}.} We use the strict models introduced in our last section. 
Let $X:=Spec\, B$ and $Y:=Spec\, C$, and $K_B$, $K_C$ the Koszul dg-algebras
of $B$ and $C$ with respect to the element $\pi$. As in our previous section, 
we have the Hopf dg-algebroid $K^2_A$ and its monoidal dg-category of
modules $\Mod^c(K_A^2)$. Now, $\Mod^c(K_A^2)$ acts on both
$\Mod^c (K_B)$ and $\Mod^c (K_C)$, the dg-categories of cofibrant $K_B$ (resp. $K_C$) dg-modules
which are perfect over $B$ (resp. over $C$). 

We define a dg-functor
$$\theta : (\Mod^c (K_B))^{o} \otimes_{\Mod^c (K^2_A)} \Mod^c (K_C) \longrightarrow 
\Mod^c (B\otimes_A C),$$
where $ \Mod^c (B\otimes_A C)$ is the category of perfect
$B\otimes_A C$-dg-modules. This dg-functor
sends a pair of objects $(E,F)$ to the
object $\mathbb{D}(E) \otimes_{K_A} F$, where $\mathbb{D}(E)=\underline{Hom}_{K_B}(E,K_B)$
is the $K_B$-linear dual of $E$. After localization with respect to 
quasi-isomorphisms, we get a well defined morphism in $\dgcat_A$
$$\theta : \Coh(X_s)^o \otimes_{\B^+} \Coh(Y_s) \longrightarrow \Coh(Z).$$

In order to finish the proof, we have to check the following two 
conditions:
\begin{enumerate}

\item The image of $\theta$ generates (by shifts, sums, cones and retracts) 
the full sub-dg-category $\Coh_{Z_s}(Z)$;

\item The dg-functor $\theta$ above is fully faithful. 

\end{enumerate}

Now, on the level of objects the dg-functor $\theta$ sends 
a pair $(E,F)$, of coherent complexes on $X_s$ and $Y_s$, to the coherent complex on $Z$
$$j_*(\mathbb{D}(E)\boxtimes_k F),$$
where $j : Z_s \hookrightarrow Z$ is the closed embedding, 
$\boxtimes_k $ denotes the external product on $Z_s=X_s\times_s Y_s$, and $\mathbb{D}(E)$ denotes the Grothendieck dual  of $E$ on the Gorenstein scheme $X_s$.
It is known that 
coherent complexes of the form 
$E\boxtimes_k F$ generate $\Coh(Z_s)$. As the coherent complexes 
of the form $j_*(G)$, for $G \in \Coh(Z_s)$, clearly generate 
the dg-category $\Coh_{Z_s}(Z)$, we see that condition $(1)$ above is satisfied.

It now remains to show that $\theta$ is fully faithful. Given two pairs 
of objects $(E,F), (E',F') \in \Coh(X_s) \times \Coh(Y_s)$, we have the morphism
induced by $\theta$
$$\rh(E',E) \otimes_{\B^+(1)} \rh(F,F') \longrightarrow
\rh(j_*(\mathbb{D}(E)\boxtimes F),j_*(\mathbb{D}(E')\boxtimes F')),$$
where $\B^+(1)$ denotes the algebra of endomorphism of the unit
in $\B^+$. As we have already seen, $\B^+\simeq k[u]$ as an $\E_1$-algebra. 
As $X_s$ and $Y_s$ are Gorenstein schemes, the structure sheaf 
$\OO$ is a dualizing complex, and $E \mapsto \mathbb{D}(E)$
is an (anti)equivalence. We thus have $\rh(E',E) \simeq \rh(\mathbb{D}(E),\mathbb{D}(E'))$,
and the above morphism can thus be written in the form
$$\rh(E,E') \otimes_{k[u]} \rh(F,F') \longrightarrow
\rh(j_*(E\boxtimes F),j_*(E'\boxtimes F')),$$
where it is simply induced by the direct image functor $j_*$.
Both the source and the target of the above morphism enter in a distinguished triangle. 
On the left hand side, for any two $k[u]$-dg-modules $M$ and $N$, we 
have a triangle of $A$-dg-modules
$$\xymatrix{
(M \otimes_k N)[-2] \ar[r] & M \otimes_k N \ar[r] & M\otimes_{k[u]} N,}$$
where the morphism on the right is the natural projection. This exact triangle
follows from the isomorphism $M\otimes_{k[u]} N \simeq (M\otimes_k N)\otimes_{k[u]\otimes_k k[u]}k[u]$, and from the exact triangle of $k[u]$ bi-dg-modules
$$\xymatrix{
(k[u]\otimes_k k[u]) [-2] \ar[r]^-{a} & k[u]\otimes_k k[u] \ar[r]^-{m} & k[u],}$$
where $a$ is given by multiplication by $(u\otimes 1 - 1 \otimes u)$, and $m$ is the product map.\\
On the right hand side, we have, by adjunction,
$$\rh(j_*(E\boxtimes F),j_*(E'\boxtimes F')) \simeq
\rh(j^*j_*(E\boxtimes F),E'\boxtimes F').$$
The adjunction map $j^*j_* \rightarrow id$, provides an exact triangle 
of coherent sheaves on $Z_s$
$$\xymatrix{
E\boxtimes_k F [1] \ar[r] & j^*j_*(E\boxtimes_k F) \ar[r] & E\boxtimes_k F.}$$
The coboundary map of this triangle  
$$E\boxtimes_k F \longrightarrow E\boxtimes_k F [2]$$
is precisely given by the action of $k[u]$. We thus obtain another exact triangle
$$\xymatrix{
(\rh(E,F)\otimes_k \rh(E',F'))[-2]  \ar[r] &  
\rh(E,F) \otimes_k \rh(E',F') \ar[r] &}$$

$$\xymatrix{ \ar[r]& \rh(j_*(E\boxtimes F),j_*(E'\boxtimes 
F')).}
$$
By inspection, the morphism $\theta$ is compatible with these two triangles, and provides
an equivalence 
$$\rh(E,E') \otimes_{k[u]} \rh(F,F') \longrightarrow
\rh(j_*(E\boxtimes F),j_*(E'\boxtimes F')).$$
Therefore $\theta$ is fully faithful.
 \hfill $\Diamond$ \\

Since $X/S$ and $Y/S$ are generically smooth, also $Z=X\times_S Y$ is generically smooth 
over $S$. Therefore $\Sing(Z) \simeq \Coh_{Z_s}(Z)/\Perf_{Z_s}(Z)$. Hence, in order to 
finish the proof of Proposition \ref{dg-kunn}, we are left to prove that the image of $
\Perf_{Z_s}(Z)$ under a quasi-inverse of  $\theta$ is generated by $\Coh(X_s)^o 
\otimes_{\B^+} \Perf(Y_s)$ and $\Perf(X_s)^o \otimes_{\B^+} \Coh(Y_s)$. 
Now, the dg-category $\Coh(X_s)^o 
\otimes_{\B^+} \Perf(Y_s)$ is generated by the image of the natural dg-functor
$$\Coh(X_s)^o 
\otimes_{A} \Perf(Y_s) \longrightarrow \Coh(X_s)^o 
\otimes_{\B^+} \Perf(Y_s),$$
and the dg-category $\Perf(X_s)^o \otimes_{\B^+} \Coh(Y_s)$ is generated by the image of the natural dg-functor
$$\Perf(X_s)^o \otimes_{A} \Coh(Y_s) \longrightarrow \Perf(X_s)^o \otimes_{\B^+} \Coh(Y_s).$$
Therefore, what we have to prove is that $\Perf_{Z_s}(Z)$ is generated by 
objects of the form $E \boxtimes_A F$, for $E \in \Coh(X_s)$, 
$F \in \Coh(Y_s)$, one of them being perfect. 

We first notice that, if $E$ or $F$ is perfect, indeed, $E\boxtimes_A F$ 
is a perfect complex on $Z$. To see this, we may localize on $Z$, and write
$X=Spec\, B$ and $Y=Spec\, C$. The objects $E$ and $F$ then correspond to bounded
coherent complexes over $B_s=B\otimes_A k$ and $C_s=C\otimes_A k$ respectively. 
Assume that $E$ is perfect (the complementary case being totally analogous). By d\'evissage we can assume that $E=B_s$, so that 
$E\boxtimes_A F \simeq (B\boxtimes_A F)\otimes_A k$. In other words, 
if we write $j : Y_s \hookrightarrow Y$ and 
$i : Z_s \hookrightarrow Z$ for the closed embeddings, and $p : Z \longrightarrow Y$
for the second projection, we have
$$E\boxtimes_A F \simeq i_*i^*p^*(j_*(F)).$$
But $j_*(F)$ is perfect on $Y$ (because $Y$ is regular), so $p^*j_* (F)$
is perfect on $Z$. Moreover, since $i$ is an lci morphism, $i_*i^*$ preserves
perfect complexes on $Z$, and we conclude that, indeed, $E\boxtimes_A F$ 
is a perfect complex on $Z$.

Finally, the fact that the image of $\Coh(X_s)^o 
\otimes_{\B^+} \Perf(Y_s)$ and $\Perf(X_s)^o \otimes_{\B^+} \Coh(Y_s)$ inside
$\Coh_{Z_s}(Z)$
generates the whole dg-category $\Perf_{Z_s}(Z)$
follows from the fact that the image of $\Perf(X_s)^o \otimes_{A} \Perf(Y_s)$
inside $\Perf(Z)$ already generates $\Perf_{Z_s}(Z)$. Indeed, for two perfect
complexes $E$ on $X_s$ and $F$ on $Y_s$, the image of $E\boxtimes_A F$ 
in $\Perf(Z)$ is of the form $i_*(E\boxtimes_k F)\oplus i_*(E\boxtimes_k F)[1]$. 
As objects of the form $E\boxtimes_k F$ generate $\Perf(Z_s)$ we are done, and
Proposition \ref{dg-kunn} is proved.
\hfill $\Box$ \\

\subsection{Saturatedness}

As a consequence of the K\"unneth formula for dg-categories
of singularities we prove the following result. 

\begin{prop}\label{psat}
Let $X$ be a regular and flat $S$-scheme.
\begin{enumerate}

\item If $X$ is proper over $S$, then the dg-category 
$\Coh(X_s)$ is proper over $\B^+$. 

\item If $X$ is proper over $S$, then the dg-category $\Sing(X_s)$ is saturated over $\B$.

\end{enumerate}
\end{prop}

\noindent \textbf{Proof.} $(1)$ We have to show that the big morphism
$$h : \widehat{\Coh(X_s)^o \otimes_A \Coh(X_s)} \longrightarrow 
\widehat{\B^+}\simeq \widehat{\Coh(G)}.$$
is small. Here $G=s\times_S s$, and the morphism $h$ is obtained
as follows. The derived scheme $G$ is a derived groupoid acting 
on $X_s$ by means of the natural projection on the last two factors
$$\mu : G\times_s X_s=s\times_S s \times_S X \longrightarrow X_s.$$
The projection on the first and third factors provides another
morphism
$$p : G\times_s X_s \longrightarrow X_s.$$
Thus, the morphisms $p$ and $\mu$ together  
define a morphism of derived schemes
$$q : G\times_s X_s \longrightarrow X_s \times_S X_s.$$
Finally, we denote by $r : G\times_s X_s \longrightarrow G$ the natural projection.
Now, for two coherent complexes
$E$ and $F$ on $X_s$, we first form the  external Hom complex 
$\mathcal{H}om_{A}(E,F)$ which is a coherent complex on $X_s \times_S X_s$, and  
we have
$$h(E,F)\simeq r_*(q^*\mathcal{H}om_{A}(E,F)).$$
This is a quasi-coherent complex on $G$. Both $q$ and 
$r$ are local complete intersection morphisms of derived schemes, and moreover
$r$ is proper. This
implies that $q^*$ and $r_*$ preserve coherent complexes, and thus that 
$h(E,F)$ is coherent on $G$. 

$(2)$ We have $\Sing(X_s)\simeq \Coh(X_s) \otimes_{\B^+}\B$, thus $(1)$ implies that 
$\Sing(X_s)$ is proper over $\B$. To prove it is smooth we need to prove that the
coevaluation big morphism $A \longrightarrow \Sing(X_s)^o \otimes_\B \Sing(X_s)$
is a small morphism. Using our K\"unneth formula for  
dg-category of singularities (Proposition \ref{dg-kunn}) this morphism corresponds
to the data of an ind-object in $\Sing(X\times_S X)$. This object is
the structure sheaf of the diagonal $\Delta_X$ inside $X\times_S X$ which
is an object in $\Sing(X\times_S X)$. This shows that the coevaluation morphism 
is a small morphism and thus that $\Sing(X_s)$ is indeed saturated over $\B$. 
\hfill $\Box$ \\

\begin{rmk} Proposition \ref{psat} (2) remains true if we only suppose that the singular locus of $X_s$ is proper over $s= \mathrm{Spec \, k}$.
\end{rmk}

\section{Application to Bloch's conductor formula with unipotent monodromy}\label{sectionBloch}

In this Section we first recall Bloch's Conductor Conjecture and the current state of the art for it, then we prove a version of Bloch's conductor under the hypothesis that the monodromy action is unipotent, where Bloch's number is replaced by a categorical Bloch class (see Definition \ref{bnum}). We are convinced that Bloch's number always agree with the categorical Bloch class but we will not prove (nor use) this fact in this paper. We are aware of a proof of this comparison in the geometric case (i.e. when the inclusion $s \to S$ has a retraction) but we  will defer this to another paper, where we hope to give the comparison also when $S$ has mixed characteristic. 

\subsection{Bloch's Conductor Conjecture}
Our base scheme 
is a discrete valuation ring $S=\Spec\, A$, with perfect residue field $k$, and fraction field $K$. 
Let $p: X \longrightarrow S$ be proper and flat morphism of finite type, and of relative dimension 
$n$. We assume that the generic fiber $X_K$ is smooth over $K$, and that $X$ is a regular scheme. We 
write $\bar{K}$ for the separable closure of $K$ (inside a fixed algebraic closure). \\In his 1985 
paper \cite{bl}, Bloch formulated the following conductor formula conjecture which is a kind of vast 
arithmetic generalization of Gauss-Bonnet formula, where an intersection theoretic (coherent) term, 
the \emph{Bloch's number}, is conjectured to be equal to an arithmetic (\'etale) term, the 
\emph{Artin conductor}. We address the reader to  \cite{bl} and \cite{ks} for more detailed 
definitions of the various objects involved in the statement.

\begin{conj}\label{cb} \emph{\textbf{[Bloch's conductor Conjecture]}}
Under the above hypotheses on $p: X \to S$, we have an equality
$$[\Delta_X,\Delta_X]_S = \chi(X_{\bar{k}},\ell) - \chi(X_{\bar{K}},\ell) - \mathsf{Sw}(X_{\bar{K}}),$$
where $\chi(Y,\ell)$ denotes the \emph{$\ql$-adic Euler characteristic} of a variety
$Y$, for $\ell$ prime to the characteristic of $k$, $\mathsf{Sw}(X_{\bar{K}})$
is the \emph{Swan conductor} of the $\mathsf{Gal}(\bar{K}/K)$-representation $H^*(X_{\bar{K}},\Ql)$,
and $[\Delta_X,\Delta_X]_S$ is \emph{Bloch's number} of $X/S$, i.e. the degree in $\mathsf{CH}_0(k)\simeq \mathbb{Z}$ of Bloch's localised self-intersection $(\Delta_X,\Delta_X)_S \in \mathsf{CH}_{0}(X_k)$ of the
diagonal in $X$. The (negative of the) rhs is called the \emph{Artin conductor} of $X/S$, and denoted by $\mathsf{Art}(X/S)$.
\end{conj}

Conjecture \ref{cb} is known to hold in several special cases that we recall below:

\begin{enumerate}

\item When $k$ is of characteristic zero, Conjecture \ref{cb} follows from the
work of \cite{kap}. When furthermore $X_s$ has only isolated singularities
the conductor formula was known as the \emph{Milnor formula} stating that the dimension
of the space of vanishing cycles equals the dimension of the Jacobian ring. 

\item When $S$ is equicharacteristic, Conjecture \ref{cb} has been proved recently
in \cite{sait}, based on Beilinson's theory of singular support of $\ell$-adic 
sheaves. The special subcase of isolated singularities already appeared in 
\cite[Exp. XVI]{sga7}.

\item When $X$ is semi-stable over $S$, i.e. the reduced special fiber $(X_s)_{red} \subset 
X$ is a normal crossing divisor, Conjecture \ref{cb} was proved 
in \cite{ks}.

\end{enumerate}

In view of the above list of results, one of the major open cases is that of 
isolated singularities in \emph{mixed characteristic}, which is the conjecture
appearing in Deligne's expos\'e \cite[Exp. XVI]{sga7}. \\

It is classical and easy to see that Conjecture \ref{cb} can be reduced to the case where $S$ is excellent and strictly henselian (so that $k$ algebraically closed)\footnote{We first reduce to the complete dvr  (hence henselian and excellent) case  by proper base-change, and then we further reduce to the excellent, \emph{strictly} henselian case.}. We will thus assume from now on that \emph{$k$ is algebraically 
closed}. \\

\subsection{An analog of Bloch's Conjecture for unipotent monodromy}

We start by introducing a \emph{categorical variant} of Bloch's number $[\Delta_X.\Delta_X]_S$ defined 
in terms of dg-categories of singularities.  

The dg-category $\Sing(X_s)$ of singularities 
of the special fiber comes equipped with its canonical 
$\B$-module structure described in Proposition \ref{cotens}. As $X$ is proper
over $S$, we know by Proposition \ref{psat} that $\Sing(X_s)$ is saturated over
$\B$. We are thus entitled to take the trace (Definition \ref{nctrace}) of the identity of $\Sing(X_s)$, which 
is a morphism 
$$A \longrightarrow \mathsf{HH}(\B/A)$$ in $\dgcat_A$.
This trace morphism is then, by definition, determined by a perfect
$\mathsf{HH}(\B /A)^o$-dg-module $\mathsf{HH}(\Sing(X_s)/\B)$, and thus provides us with a class in $K$-theory
$$[\mathsf{HH}(\Sing(X_s)/\B)] \in K_0(\mathsf{HH}(\B /A)).$$
The Chern character natural transformation (applied to $\mathsf{HH}(\B /A)$, see Definition \ref{d3}) of this class is an element $Ch_{\ell, \mathsf{HH}(\B /A)}([\mathsf{HH}(\Sing(X_s)/\B)])$ in $\pi_0 |r_{\ell}(\mathsf{HH}(\B /A))|= H^0(S_{\textrm{\'et}},r_{\ell}(\mathsf{HH}(\B /A)))$.

\begin{df}\label{bnum}
With the above notations and hypotheses on $X/S$, the 
\emph{categorical Bloch class} of $X/S$ is defined as
$$[\Delta_X,\Delta_X]_S^{cat}:=Ch_{\ell}([\mathsf{HH}(\Sing(X_s) / \B]) \in H^0(S_{\textrm{\'et}}, r_{\ell}(\mathsf{HH}(\B/A))).$$
\end{df}

Notice that though $\B$ is just an associative algebra in $\dgcat_A$ (an $E_2$-algebra over $A$), its $\ell$-adic realization  $r_{\ell}(\B)$ is in fact a commutative monoid in $Sh_{\Ql}(S)$, naturally equivalent to 
$i_*(\Ql(\beta)\oplus \Ql(\beta)[1])$ where $i: s \to S$ is the inclusion of the closed point.
Therefore there is a canonical augmentation $\mathsf{a}: \mathsf{HH}(r_{\ell}(\B)/r_{\ell}(A)) \to r_{\ell}(\B)$, hence an induced augmentation $$\mathsf{a}: H^0(S_{\textrm{\'et}}, \mathsf{HH}(r_{\ell}(\B)/r_{\ell}(A)))=\pi_0|\mathsf{HH}(r_{\ell}(\B)/r_{\ell}(A))| \longrightarrow
 \pi_0|r_{\ell}(\B)|=H^0(S_{\textrm{\'et}},r_{\ell}(B)) \simeq \Ql$$ that is a left inverse to $$\sigma: \Ql \simeq H^0(S_{\textrm{\'et}},r_{\ell}(\B)) \to H^0(S_{\textrm{\'et}}, \mathsf{HH}(r_{\ell}(\B)/r_{\ell}(A)))$$ (induced by $r_{\ell}(\B) \to \mathsf{HH}(r_{\ell}(\B)/r_{\ell}(A))$).\\
 
For $\lambda \in \Ql$, we will simply write $\lambda^{\wedge}$ for the image of $\lambda$ via the composition $$\xymatrix{\Ql \simeq H^0(S_{\textrm{\'et}},r_{\ell}(\B)) \ar@{^{(}->}[r]^-{\sigma} & H^0(S_{\textrm{\'et}}, \mathsf{HH}(r_{\ell}(\B)/r_{\ell}(A))) \ar[r]^-{\alpha} & H^0(S_{\textrm{\'et}}, r_{\ell}(\mathsf{HH}(\B/A))).}$$

We are now ready to prove our version of Bloch's Conductor Conjecture for unipotent monodromy..

\begin{thm}\label{blochthm} Let $X/S$ be as in Conjecture \ref{cb}, and assume further that 
the inertia subgroup $\mathrm{I}:=\mathsf{Gal}(\bar{K}/K^{\textrm{unr}}) \subseteq 
\mathsf{Gal}(\bar{K}/K) $ acts unipotently on $H^*(X_{\bar{K}}, \ql)$. Then we have an equality 
$$[\Delta_X,\Delta_X]_S^{cat}=\chi(X_{\bar{k}},\ell)^{\wedge} - \chi(X_{\bar{K}},\ell)^{\wedge}$$ in $H^0(S_{\textrm{\'et}}, r_{\ell}(\mathsf{HH}(\B/A)))$.
\end{thm}

\begin{rmk}\label{tomakeclear}
Since we have not proved that the map $$\alpha: H^0(S_{\textrm{\'et}}, \mathsf{HH}(r_{\ell}(\B)/r_{\ell}(A))) \to H^0(S_{\textrm{\'et}}, r_{\ell}(\mathsf{HH}(\B/A)))$$ is injective, the equality in Theorem \ref{blochthm} is not a priori an equality of integers (and it might even be the trivial equality $0=0$ of classes in $H^0(S_{\textrm{\'et}}, r_{\ell}(\mathsf{HH}(\B/A)))$). However, we conjecture that the canonical map $u: \B \to \mathsf{HH}(\B / A)$ is in fact an $\mathbb{A}^1$-homotopy equivalence; this would imply that $r_{\ell}(u): r_{\ell}(\B) \to r_{\ell}(\mathsf{HH}(\B / A))$ is an equivalence in $Sh_{\Ql}(S)$, so that $H^0(S_{\textrm{\'et}}, r_{\ell}(\mathsf{HH}(\B/A))) \simeq H^0(S_{\textrm{\'et}}, r_{\ell}(\B)) \simeq \Ql$. In particular $\alpha$ would be injective, and the equality in Theorem \ref{blochthm} would indeed be an equality of integers.

\end{rmk}

\begin{rmk}\label{unipimpliestame} Recall that the inertia group $\mathrm{I}$ sits in an extension of pro-finite
groups
$$\xymatrix{
1 \ar[r] & P \ar[r] & \mathrm{I} \ar[r] & \mathrm{I}_t \ar[r] & 1\, ,}$$
where $\mathrm{I}_t$ is the \emph{tame inertia} quotient, and the \emph{wild inertia} $P$ is a pro-$p$-subgroup. For any continuous finite dimensional $\Ql$-representation
$V$ of $\mathrm{I}$, the group $P$ acts through a finite quotient on $V$ (see \cite{}). Therefore, if $\mathrm{I}$ is supposed to act unipotently on $V$, then $P$ acts necessarily trivially, i.e. the $\mathrm{I}$-action factors through a $\mathrm{I}_t$-action, i.e., by definition, the $\mathrm{I}$ action is \emph{tame}. By definition, the Swan conductor $Sw_{\mathrm{I}}(V)$(see e.g. \cite[\S 6.1]{ks}) vanishes for a tame $\mathrm{I}$-representation  $V$. As a consequence, granting Remark \ref{tomakeclear}, under the hypothesis of unipotent action of $\mathrm{I}$ on $H^*(X_{\bar{K}}, \ql)$, Conjecture \ref{cb} becomes Theorem \ref{blochthm} \emph{if} Bloch's number $[\Delta_X,\Delta_X]_S$ is replaced by the categorical Bloch class $[\Delta_X,\Delta_X]_S^{cat}$. 
Though we currently know a proof only in the geometric case,  we are actually convinced that 
the categorical Bloch class is  \emph{always}, i.e. regardless any hypothesis on the monodromy action and on the mixed or equal characteristic property of $S$, an integer equal to Bloch's intersection number $[\Delta_X,\Delta_X]_S$
appearing in Conjecture \ref{cb}. This fact will not be used in this paper and will be more closely investigated in a future one.
\end{rmk}

\noindent \textbf{Proof of Thm. \ref{blochthm}.} We may suppose $S$ strictly henselian, so that $\mathrm{I}= \mathsf{Gal}(\bar{K}/K)$. We want to apply our trace formula (Theorem \ref{ttrace}) to $T=\Sing(X_s)$ and 
$f=id$. In order to do this, we need to check that the conditions of being \emph{saturated} and \emph{$\ell^{\otimes}$-admissible} over $\B$ are met by such $T$.

By Proposition  \ref{psat}, we know that $T$ is \emph{saturated over} $\B$.\\
Let us now show that $T$ is also \emph{$\ell^{\otimes}$-admissible} over $\B$, i.e. that the canonical map 
$$\varphi: r_{\ell}(T^o)\otimes_{r_{\ell}(\B)}r_{\ell}(T) \longrightarrow
r_{\ell}(T^o\otimes_{\B}T)$$ is an equivalence.\\
First of all we notice that, since both the source and the target of $\varphi$ 
are $\ell$-adic complexes on $S$, supported at $s$, it is enough to show that
$i^*(\varphi)$ is an equivalence, $i: s \to S$ being the canonical inclusion (note that $i^* i_*$ is equivalent to the identity functor).\\
Let us first elaborate on the target of $i^*(\varphi)$. By K\"unneth for dg-categories of singularities, we have the canonical
equivalence of Proposition \ref{dg-kunn}
$$T^o \otimes_\B T \simeq \Sing(X\times_S X).$$
Moreover, since a unipotent 
action of $\mathrm{I}$ is also tame, by Cor. \ref{ckunn} we have 
$$r_{\ell}(\Sing(X\times_S X)) \simeq 
(p_{Z_s})_{*}(\nu_X \boxtimes \nu_X)^{\mathrm{I}} (\beta),$$
where $p_{Z_s}: Z_s=X_s \times_{s} X_s \to S$ is the composite $\xymatrix{Z_s \ar[r]^-{q_{Z_s}} & s \ar[r]^-{i} & S}$, and, as usual, $\nu_X$ denotes vanishing cycles for $X/S$ on $X_s$. Therefore 
$$i^*r_{\ell}(\Sing(X\times_S X)) \simeq 
(q_{Z_s})_{*}(\nu_X \boxtimes \nu_X)^{\mathrm{I}} (\beta)$$
in $Sh_{\Ql}(s)$. Since $q_{Z_s}= q_{X_s} \times_s q_{X_s}$, where $q_{X_s}: X_s \to s$ is the canonical map, we have 
\begin{equation} 
i^*r_{\ell}(\Sing(X\times_S X)) \simeq 
(q_{Z_s})_{*}(\nu_X \boxtimes \nu_X)^{\mathrm{I}} (\beta) \simeq ((q_{X_s})_{*}(\nu_X) \otimes_{\Ql}(q_{X_s})_{*}(\nu_X))^{\mathrm{I}}(\beta) 
\end{equation}
in $Sh_{\Ql}(s)$, that we may rewrite as
\begin{equation} \label{RHS}
i^*r_{\ell}(\Sing(X\times_S X)) \simeq (\mathbb{H}(X_s, \nu_X) \otimes_{\Ql}\mathbb{H}(X_s, \nu_X))^{\mathrm{I}}(\beta) 
\end{equation}
where we have introduced the hypercohomology $\mathbb{H}(X_s, \mathcal{E}):= (q_{X_s})_{*}(\mathcal{E})$, for $\mathcal{E} \in Sh_{\Ql}(X_s)$.\\
Let us now look more carefully at the source of the map $i^*(\varphi)$. 
By the main theorem of \cite{brtv} (see also formula (\ref{prelim})), we have 
$$r_{\ell}(T) \simeq (p_{X_s})_{*}(\nu_X[-1]^{\mathrm{I}})(\beta)$$
where $p_{X_s}: X_s \to S$ is the composite $\xymatrix{X_s \ar[r]^-{q_{X_s}} & s \ar[r]^-{i} & S}$. Therefore,
\begin{equation}
i^*(r_{\ell}(T^o)\otimes_{r_{\ell}(\B)}r_{\ell}(T)) \simeq i^*(r_{\ell}(T^o))\otimes_{i^*r_{\ell}(\B)}i^*r_{\ell}(T) \simeq ((q_{X_s})_{*}(\nu_X[-1])^{\mathrm{I}} \otimes_{\Ql^{\mathrm{I}}} q_{X_s})_{*}(\nu_X[-1])^{\mathrm{I}})(\beta)
\end{equation} 
in $Sh_{\Ql}(s)$, that we may rewrite as\footnote{Note that $r_{\ell}(T) \simeq r_{\ell}(T^o)$ because the $K$-theory of a dg-category is canonically isomorphic to the $K$-theory of the opposite dg-category.}
\begin{equation}\label{LHS}
i^*(r_{\ell}(T^o)\otimes_{r_{\ell}(\B)}r_{\ell}(T)) \simeq (\mathbb{H}(X_s, \nu_X[-1])^{\mathrm{I}} \otimes_{\Ql^{\mathrm{I}}} \mathbb{H}(X_s, \nu_X[-1])^{\mathrm{I}})(\beta).
\end{equation}

By recalling that $\mathcal{E}(\beta) \simeq \mathcal{E}[-2](\beta)$, for any $\mathcal{E} \in Sh_{\Ql}(X_s)$, and by combining (\ref{RHS}) and (\ref{LHS}), we see that 
$$i^*(\varphi): i^*(r_{\ell}(T^o)\otimes_{r_{\ell}(\B)}r_{\ell}(T)) \longrightarrow
i^*r_{\ell}(T^o\otimes_{\B}T)$$
is then equivalent to $\psi[-2](\beta)$, where $\psi$ is the lax-monoidal structure morphism for the functor $(-)^{\mathrm{I}}$, applied to the pair of $\ell$-adic complexes $(\mathbb{H}(X_s, \nu_X), \mathbb{H}(X_s, \nu_X))$,
\begin{equation}
    \psi: \mathbb{H}(X_s, \nu_X)^{\mathrm{I}} \otimes_{\Ql^{\mathrm{I}}} 
\mathbb{H}(X_s, \nu_X)^{\mathrm{I}} \longrightarrow (\mathbb{H}(X_s, \nu_X) \otimes_{\Ql}
\mathbb{H}(X_s, \nu_X))^{\mathrm{I}}. 
\end{equation}

The fact that the morphism $\psi$ is an equivalence is a consequence
of the following general lemma.

\begin{lem}\label{luni}
Let $\D_{uni}(\mathrm{I},\Ql)$ be the full sub-$\s$-category of 
$\D_c(Spec\, K,\Ql) \simeq \D^{\mathrm{I}}_c(Spec\, \bar{K},\Ql) $ consisting of all objects $E$ for which
the action of $\mathrm{I}$ on each $H^i(E)$ is unipotent. Then the invariant $\s$-functor induces
an equivalence of symmetric monoidal $\s$-categories
$$(-)^{\mathrm{I}} : \D_{uni}(\mathrm{I},\Ql) \simeq \D_{\textrm{perf}}(\Ql[\epsilon_1]),$$
where $\Ql[\epsilon_1]$ is the free commutative $\Ql$-dg-algebra generated by 
$\epsilon_1$ in degree $1$, and $\D_{\textrm{perf}}(\Ql[\epsilon_1])$ 
is its $\s$-category of perfect dg-modules.
\end{lem}
 
\noindent \textit{Proof of lemma \ref{luni}.} This is a well known fact. The $\s$-functor
$$(-)^\mathrm{I} : \D_c(Spec\, K,\Ql) \longrightarrow \D_c(\Ql)$$
is lax symmetric monoidal, so it induces a lax monoidal $\s$-functor
$$(-)^\mathrm{I} : \D_c(Spec\, K,\Ql) \longrightarrow \D(\Ql^\mathrm{I}).$$
It is easy to see that $\Ql^\mathrm{I}$ is canonically equivalent to
$\Ql[\epsilon_1]$, and the choice of such an equivalence only depends
on the choice of a generator of $H^1(\mathrm{I},\Ql)\simeq \Ql$. We thus have 
an induced lax symmetric monoidal $\s$-functor
$$(-)^\mathrm{I} : \D_c(Spec\, K,\Ql) \longrightarrow \D(\Ql[\epsilon_1]).$$
The above $\s$-functor is in fact symmetric monoidal, as it preserves unit objects,
and moreover \\$\D_c(Spec\, K,\Ql)$ is generated by the unit object $\Ql$.  \\
The lemma is then a direct consequence of the following fact: let $x \in \mathcal{C}$
be a compact object in a compactly generated stable $\s$-category $\mathcal{C}$, 
then the $\s$-functor
$$\mathcal{C}(x,-) : \mathcal{C} \longrightarrow \textsf{Mod}_{\mathbf{Sp}}(End(x))$$
induces an equivalence of stable $\s$-categories
$$\langle x \rangle  \simeq \textsf{Perf}_{\mathbf{Sp}}(End(x)),$$
where $\langle x \rangle  \subset \mathcal{C}$  denotes the full stable $\s$-category generated
by the object $x$, and $End(x)$ denotes the ring spectrum of endormophisms of $x$. 
\hfill $\diamondsuit$ \\

The above lemma implies that the map
$$\psi[-2]: \mathbb{H}(X_s,\nu_X[-1])^{\mathrm{I}} \otimes_{\Ql^{\mathrm{I}}} 
\mathbb{H}(X_s,\nu_X[-1])^{\mathrm{I}} \longrightarrow (\mathbb{H}(X_s,\nu_X[-1])\otimes_{\Ql}
\mathbb{H}(X_s,\nu_X[-1]))^{\mathrm{I}}$$
is an equivalence, and thus, as observed above, the same is true for $i^*(\varphi)$, and hence for the admissibility map
$$\varphi: r_{\ell}(T^o)\otimes_{r_{\ell}(\B)}r_{\ell}(T) \longrightarrow
r_{\ell}(T^o\otimes_{\B}T),$$ so that $T$ is indeed \emph{$\ell^{\otimes}$-admissible over} $\B$ as we wanted. \\

Now that we have checked all the conditions for the trace formula of Theorem \ref{ttrace} to hold in our case, 
we get an equality
\begin{equation}\label{now}
Ch_{\ell}([\mathsf{HH}(T/\B, id)]) = tr_{\rl(\B)}(id: r_{\ell}(T) \to r_{\ell}(T))
\end{equation} 
in $H^0(S_{\textrm{\'et}},r_{\ell}(\mathsf{HH}(\B /A)))$. Recall that $$tr_{\rl(\B)}(id: r_{\ell}(T) \to r_{\ell}(T)):= \alpha(Tr_{\rl(\B)}(id: r_{\ell}(T) \to r_{\ell}(T))),$$ where $\alpha: H^0(S_{\textrm{\'et}}, \mathsf{HH}(r_{\ell}(\B)/r_{\ell}(A))) \to H^0(S_{\textrm{\'et}}, r_{\ell}(\mathsf{HH}(\B/A)))$ is the canonical map.  \\

We are left to 
identify the two sides of (\ref{now}).\\The left hand side is,
by definition, our categorical Bloch number $[\Delta_X,\Delta_X]_S^{cat}$. 

Let us now investigate the r.h.s. of (\ref{now}): we need to show that  $Tr_{r_{\ell}(B)}(r_{\ell}(id_T))=\chi(X_{\bar{k}},\ell) - \chi(X_{\bar{K}},\ell)$. \\ As already noticed, $r_{\ell}(\B)$ is a commutative monoid in $Sh_{\Ql}(S)$, naturally equivalent to 
$i_*(\Ql(\beta)\oplus \Ql(\beta)[1])$ where $i: s \to S$ is the inclusion of the closed point.
Therefore there is a canonical augmentation $\mathsf{a}: \mathsf{HH}(r_{\ell}(\B)/r_{\ell}(A)) \to r_{\ell}(\B)$, hence an induced augmentation $$\mathsf{a}: H^0(S_{\textrm{\'et}}, \mathsf{HH}(r_{\ell}(\B)/r_{\ell}(A)))=\pi_0|\mathsf{HH}(r_{\ell}(\B)/r_{\ell}(A))| \longrightarrow
 \pi_0|r_{\ell}(\B)|=H^0(S_{\textrm{\'et}},r_{\ell}(B)).$$

Consider the diagram  $$\xymatrix{H^0(S_{\textrm{\'et}}, \mathsf{HH}(r_{\ell}(\B)/r_{\ell}(A))) \ar[r]^-{\alpha} \ar@/_0.4cm/[d]_-{\mathsf{a}} & H^0(S_{\textrm{\'et}}, r_{\ell}(\mathsf{HH}(\B/A)))\\ \Ql \simeq H^0(S_{\textrm{\'et}}, r_{\ell}(\B)) \ar[u]_-{\sigma}  \ar[ru]_-{\textrm{can}} &  }$$ where can is the map induced by the canonical map $u:\B \to \mathsf{HH}(\B/A)$, and $\sigma$ is the map induced by $r_{\ell}(u)$, so that $\mathsf{a}\circ \sigma= \mathrm{id}$.
By compatibility between the non-commutative trace and the commutative trace (Remark \ref{tr-commvsnoncomm}) $$\mathsf{a}(Tr_{r_{\ell}(\B)}(r_{\ell}(id_T))= Tr^{\textrm{c }}_{r_{\ell}(\B)}(r_{\ell}(id_T))$$.\\
Since $r_{\ell}(\B)\simeq i_*(\Ql(\beta)\oplus \Ql(\beta)[1])$, we have $K_0(r_{\ell}(\B)) \simeq \mathbb{Z}$ and the following commutative diagram $$\xymatrix{\Ql \simeq H^0(S_{\textrm{\'et}},r_{\ell}(\B)) \ar@{^{(}->}[r]^-{\sigma} & H^0(S_{\textrm{\'et}}, \mathsf{HH}(r_{\ell}(\B)/r_{\ell}(A))) \ar[r]^-{\mathsf{a}} &  H^0(S_{\textrm{\'et}},r_{\ell}(\B)) \\
K_0(r_{\ell}(\B)) \ar[ur]^-{Tr (id_{-}) } \ar[urr]^-{Tr^{\textrm{c}} (id_{-})} & & \\
 & \mathbb{Z} \ar[ul]^{\textrm{iso}} \ar@/^0.4cm/[uul] \ar@/_1cm/[uu] \ar[uur] &  }
$$ where the maps with source $\mathbb{Z}$ are the unique maps of commutative rings, and $\mathsf{a} \circ Tr(id_{-})= Tr^{\textrm{c}}(id_{-})$ again by the functorial compatibility of non-commutative and commutative traces. By considering the class $[r_{\ell}(T)]$ in $K_0(r_{\ell}(\B))$, and using that $\mathsf{a}\circ \sigma= \mathrm{id}$, we deduce from the previous commutative diagram the equality of $\ell$-adic numbers
\begin{equation}\label{sigma}\sigma(Tr^{\textrm{c}}_{r_{\ell}(\B)}(r_{\ell}(id_T))= Tr_{r_{\ell}(B)}(r_{\ell}(id_T)).\end{equation}

Now, unfolding the definition of $Tr_{\rl(\B)}$, 
and using Lemma \ref{luni}, the $\ell$-adic number (\ref{sigma}) can be described as follows. 
The dg-algebra $\Ql^\mathrm{I}$ is such that $K_0(\Ql^\mathrm{I})\simeq \mathbb{Z}$. Viewing
$\mathbb{Z}$ inside $\Ql$, this isomorphism is induced by 
sending the class of dg-module $E$ to the trace of the identity
inside $\mathsf{HH}_0(\Ql^\mathrm{I})\simeq \Ql$. Using the functoriality of 
the trace for the morphism of commutative dg-algebras $\Ql^\mathrm{I} \longrightarrow
\Ql^\mathrm{I}(\beta)$, we see that (\ref{sigma}) is simply the trace
of the identity on $\mathbb{H}(X_s,\nu_X[-1])$, as an object
in $\D_{uni}(\mathrm{I},\Ql)$. This trace is easy to compute, 
as it equals $1$ on the unit object $\Ql$. Since the unit object
generates $\D_{uni}(\mathrm{I},\Ql)$, we have that the trace of the identity on 
any object $E \in \D_{uni}(\mathrm{I},\Ql)$ equals the Euler characteristic of the underlying complex
of $\Ql$-vector spaces. 

We thus have shown that (\ref{sigma}) can be rewritten as
$$Tr_{r_{\ell}(B)}(r_{\ell}(id_T))= \sum_{i}(-1)^i \dim_{\Ql}\, H^{i-1}(X_s,\nu_X).$$
However, by proper base change, the complex $\mathbb{H}(X_s,\nu_X)$
appears in an exact triangle
$$\xymatrix{
H(X_{\bar{k}},\Ql) \ar[r] & H(X_{\bar{K}},\Ql) \ar[r] & 
\mathbb{H}(X_s,\nu_X)}$$
and thus we have the equality
$$Tr_{r_{\ell}(B)}(r_{\ell}(id_T))=\chi(X_{\bar{k}},\ell) - \chi(X_{\bar{K}},\ell)$$
in $\mathbb{Z} \hookrightarrow \Ql \hookrightarrow H^0(S_{\textrm{\'et}}, \mathsf{HH}(r_{\ell}(\B)/r_{\ell}(A)))$.\\
Therefore 

$$[\Delta_X,\Delta_X]_S^{cat}=\chi(X_{\bar{k}},\ell)^{\wedge} - \chi(X_{\bar{K}},\ell)^{\wedge},$$
as required. \hfill $\Box$ \\

\begin{appendix}

\section{Localizations of monoidal dg-categories}

In this Appendix we remind some basic facts about localizations of dg-categories
as introduced in \cite{to-dgcat}. The purpose of the section is to explain the 
multiplicative
properties of the localization construction. In particular, we explain how
localization of strict monoidal dg-categories gives rise to monoids in $\dgcat_A$, 
and thus to monoidal dg-categories in the sense of our Definition \ref{weakmon}.\\

Let $T$ be a dg-category over $A$, together with $W$ a set of morphisms
in $Z^0(T)$, the underlying category of $T$ (this is the category of $0$-cycles
in $T$, i.e. $Z^0(T)(x,y):= Z^0 (T(x,y))$). For the sake of brevity, we will just say 
that \emph{$W$ is a set of maps in $T$}. In other words, we allow $W$ not to be
strictly speaking a \emph{subset} of the set of morphisms in $T$, but 
just a set together with a map $W \rightarrow Mor(T)$ from $W$ to the set 
of morphisms in $T$.
Recall that a localization
of $T$ with respect to $W$, is a dg-category $L_WT$ together with a morphism 
in $\dgcat_A$
$$l : T \longrightarrow  L_WT$$
such that, for any $U \in \dgcat_A$, map induced by $l$ on mapping spaces
$$Map(L_WT,U) \longrightarrow Map(T,U)$$
is fully faithful and its image consists of all $T \rightarrow U$ sending
$W$ to equivalences in $U$ (i.e. the induced functor $[T] \rightarrow [U]$
sends elements of $W$ to isomorphisms in $[U]$).

As explained in \cite{to-dgcat}, localization always exists, and are unique up to a contractible
space of choices (because they represents an obvious $\s$-functor). We will describe here
a model for $(T,W) \mapsto L_WT$ which will have nice properties with respect to
tensor products of dg-categories. For this, let $dgcat_A^{W,c}$ be the category of
pairs $(T,W)$, where $T$ is a dg-category with cofibrant hom's over $A$, and
$W$ a set of maps in $T$. Morphisms $(T,W) \longrightarrow (T',W')$
in $dgcat_A^{W,c}$ are dg-functors $T \longrightarrow T'$ sending $W$ to $W'$.

We fix once for all a factorization
$$\xymatrix{
\Delta^1_A \ar[r]^-j & \tilde{I} \ar[r]^-{p} & \overline{\Delta}^1_A,}$$
with $j$ a cofibration and $p$ a trivial fibration. Here $\Delta^1_A$ is the
$A$-linearisation of the category $\Delta^1$ that classifies morphisms, and
$\overline{\Delta}^1_A$ is the linearisation of the category that classifies
isomorphisms. For an object $(T,W) \in dgcat^{W,c}_A$ we define
$W^{-1}T$ by the following cocartesian diagram in dg-categories
$$\xymatrix{
\coprod_W \Delta^1_A \ar[r] \ar[d] & T \ar[d] \\
\coprod_W \tilde{I} \ar[r] & W^{-1}T,}$$
where $\coprod_W \Delta^1_A \longrightarrow T$ is the canonical dg-functor
corresponding to the set $W$ of morphisms in $T$. 

\begin{lem}
The canonical morphism $l : T \longrightarrow W^{-1}T$ defined above is
a localization of $T$ along $W$.
\end{lem}

\noindent \textbf{Proof.} According to \cite{to-dgcat}, the localization of $T$ can 
be construced as the homotopy push-out of dg-categories
$$\xymatrix{
\coprod_W \Delta^1_A \ar[r] \ar[d] & T \ar[d] \\
\coprod_W A \ar[r] & L_WT.}$$
The lemma then follows from the observation that when $T$ has cofibrant hom's, then 
the push-out diagram defining $W^{-1}T$ is in fact a homotopy push-out diagram. 
\hfill $\Box$ \\

The construction $(T,W) \longrightarrow W^{-1}T$ clearly defines
a functor
$$dgcat^{W,c}_A \longrightarrow dgcat_A^c$$
from $dgcat^{W,c}_A$ to $dgcat_A^c$, the category of dg-categories with cofibrant hom's. 
Moreover, this functor comes equipped with a natural symmetric colax monoidal structure. Indeed, 
$dgcat^{W,c}_A$ is a symmetric monoial category, where the tensor product is given by 
$$(T,W) \otimes (T',W'):=(T\otimes_A T',W\otimes id \cup id \otimes W').$$
We have a natural map
$$T\otimes_A T' \longrightarrow (W^{-1}T) \otimes_A ((W')^{-1}T'),$$
which by construction has a canonical extension
$$(W\otimes id \cup id \otimes W')^{-1}(T\otimes_A T') \longrightarrow (W^{-1}T) \otimes_A
((W')^{-1}T').$$
The unit in $dgcat_A^{W,c}$ is $(A,\emptyset)$,
which provides an canonical isomorphism $(\emptyset)^{-1}A\simeq A$. These data 
endow the functor $(T,W) \mapsto W^{-1}T$ with a symmetric colax monoidal structure. 
By composing with the canonical symmetric monoidal $\s$-functor
$dgcat_A^c \longrightarrow \dgcat_A,$
we get a symmetric colax monoidal $\s$-functor
$$dgcat_A^{W,c} \longrightarrow \dgcat_A,$$
which sends $(T,W)$ to $W^{-1}T$. By \cite[Ex. 4.3.3]{dgcat}, this colax symmetric monoidal 
$\s$-functor is in fact monoidal. We thus have a symmetric monoidal localization $\s$-functor 
$$W^{-1}(-) : dgcat_A^{W,c} \longrightarrow \dgcat_A.$$

As a result, if $T$ is a (strict) monoid in $dgcat_A^{W,c}$, then 
$W^{-1}T$ carries a canonical structure of a monoid in $\dgcat_A$. This applies particularly
to stict monoidal dg-categories endowed with a compatible notion of equivalences. By MacLane coherence
theorem, any such a structure can be turned into a strict monoid in $dgcat_A^{W,c}$, and by localization
into a monoid in $\dgcat_A$. In other words, the localization of 
a monoidal dg-categorie along a set of maps $W$ that is compatible with the monoidal
structure, is a monoid in $\dgcat_A$. The same is true for dg-categories which are
modules over a given monoidal dg-category.

\end{appendix}

\bigskip
\bigskip

\noindent
Bertrand To\"{e}n, {\sc IMT, CNRS, Universit\'e de Toulouse},
Bertrand.Toen@math.univ-toulouse.fr 

\smallskip

\noindent
Gabriele Vezzosi, {\sc DIMAI, Universit\`a di Firenze},
gabriele.vezzosi@unifi.it

\end{document}